\documentclass[11pt,a4paper]{article}

\usepackage[utf8]{inputenc}
\usepackage{amsmath}
\usepackage{amsfonts}
\usepackage[colorlinks,linkcolor=blue,anchorcolor=blue,citecolor=blue]{hyperref}
\usepackage{amssymb}
\usepackage{extarrows}
\usepackage{algorithmic}
\usepackage{algorithm}
\usepackage{caption}
\usepackage{amsthm}
\usepackage{comment}
\usepackage{bm}
\usepackage{tikz}
\usetikzlibrary{arrows.meta}       
\usetikzlibrary{positioning}       
\usetikzlibrary{calc}              
\usetikzlibrary{shapes}  
\usepackage{verbatim}
\usepackage{hyperref}
\usepackage{cite}
\usepackage[numbers]{natbib}
\usepackage{enumerate}
\PassOptionsToPackage{normalem}{ulem}
\usepackage{ulem}
\usepackage[margin=1in]{geometry} 
\usepackage{array}
\usepackage{cases}
\usepackage{mathrsfs}
\usepackage{longtable}
\allowdisplaybreaks[4]
\allowdisplaybreaks
\newcolumntype{L}[1]{>{\raggedright\let\newline\\\arraybackslash\hspace{0pt}}m{#1}}
\newcolumntype{C}[1]{>{\centering\let\newline\\\arraybackslash\hspace{0pt}}m{#1}}
\newcolumntype{R}[1]{>{\raggedleft\let\newline\\\arraybackslash\hspace{0pt}}m{#1}}

\usepackage{color}
\usepackage{dsfont}
\usepackage{marginnote}
\usepackage{enumitem}
\usepackage{graphicx}
\usepackage{comment}
\usepackage{float}
\usepackage{subfigure}

\theoremstyle{plain}
\newtheorem{theorem}{\protect Theorem}[section]
\newtheorem{prop}[theorem]{\protect Proposition}

\newtheorem{lemma}[theorem]{\protect Lemma}
\newtheorem{remark}[theorem]{\protect Remark}
\newtheorem{ass}{\protect Assumption}
\newtheorem{corollary}[theorem]{\protect Corollary}

\makeatletter
\newenvironment{breakablealgorithm}[1]
{
	\par\addvspace{6pt}
	\refstepcounter{algorithm}
	\addcontentsline{loa}{algorithm}
	{\protect\numberline{\thealgorithm}#1}
	\noindent
	\hrule height 0.8pt
	\vspace{3pt}
	\noindent\textbf{Algorithm \thealgorithm: #1}\par
	\vspace{3pt}
	\hrule height 0.4pt
	\vspace{3pt}
}
{
	\par
	\vspace{3pt}
	\noindent
	\hrule height 0.8pt
	\par\addvspace{6pt}
}
\makeatother

\def\d{\mathrm{d}}

\newcommand{\R}{\mathbb{R}}
\newcommand{\E}{\mathbb{E}}
\newcommand{\Sgn}{{\rm sgn}}

\newcommand{\U}{\mathscr{U}}

\newcommand{\F}{\mathcal{F}}
\newcommand{\Fb}{\mathbb{F}}
\newcommand{\Pb}{\mathbb{P}}

\newcommand{\Law}{\mathcal{L}}

\newcommand{\pa}{\partial}

\title{Model-free Reinforcement Learning for Continuous Time Optimal Stopping: A Stochastic Maximum Principle Approach}

\author{Lijun Bo\thanks{Email: lijunbo@ustc.edu.cn, School of Mathematics and Statistics, Xidian University, Xi'an, 710126, China.}
	\and
	Yijie Huang\thanks{Email: yijie.huang@polyu.edu.hk, Department of Applied Mathematics, The Hong Kong Polytechnic University, Kowloon, Hong Kong, China.}
	\and
	Jingfei Wang \thanks{Email: wjf2104296@mail.ustc.edu.cn, School of Mathematical Sciences, University of Science and Technology of China, Hefei, 230026, China.}
}
\date{\vspace{-0.3cm}}

\begin{document}
	\maketitle
	\begin{abstract}
		
		In this paper, we develop a reinforcement learning (RL) method based on stochastic maximum principle (SMP) for model-free continuous time optimal stopping problem. We consider an intensity control problem which is value-function equivalent to the original stopping problem and establish an optimal admissible intensity of the truncated intensity control problem. We derive a performance-difference formula in terms of the stopping gap for the intensity control problem using the SMP. Following the form of optimal capped intensity and the performance-difference formula, we propose an exact policy iteration for the capped intensity and prove its convergence to an optimal intensity with an explicit factorial convergence rate. We also show that the stopping gap can be learned from data by minimizing a weighted $L^2$-loss over a parameterized function. We then use the learned stopping gap to drive a relaxed policy iteration on intensity and establish its theoretical convergence. Lastly, we convert the capped learned intensity back into a near-optimal Brownian stopping time of the prime problem. Numerical experiments demonstrate the effectiveness of the proposed RL algorithm.
		
		\ \\
		
	\end{abstract}
	
	\section{Introduction}
	
	Optimal stopping is a well-studied class of decision problems in which an agent chooses a stopping time to maximize the expected reward accumulated up to that time (see, e.g., Shiryaev \cite{Shiryaev08}). In continuous time, the underlying state process is typically modeled as a diffusion that solves a stochastic differential equation. When the model coefficients of state process and reward functional are known, the dynamic programming principle (DPP) reduces optimal stopping to the solution of a variational inequality and the identification of a free boundary (see, e.g., Friedman \cite{Friedman82} and Peskir and Shiryaev \cite{PeskirShiryaev06}). The continuous time optimal stopping has been applied in financial economics such as irreversible investment (\cite{DixitPindyck08, DeAngelisETAL17}), option pricing (\cite{Glasserman04, Jacka91, Peskir05}), management science (\cite{Brancoetal12, Bochenzhang23}) and quickest detection (\cite{GloverPeskir25, PoorHadjiladis09, Shiryaev10}). However, when the full knowledge of the model primitives is unknown, reinforcement learning (RL) offers a data-driven alternative: the agent interacts with the unknown environment, observes realized states and rewards, and improves a stopping rule without estimating model coefficients.
	
	The existing continuous-time RL literature for optimal stopping has developed along a common methodological line. Dong \cite{Dong24} developed the exploratory framework as in Wang and Zhou \cite{WangZariXYZ2020}, in which the agent takes actions randomly conditioning on the current state, a regularization term is added to the reward functional, and the optimal stopping problem is reduced to a standard optimal control problem. Dai et al. \cite{Daisunxuzhou26} studied optimal stopping for diffusions with unknown primitives by penalizing the variational inequality and transforming the problem into a stochastic control problem with two actions. They randomized the control into Bernoulli distributions and added an entropy regularizer to encourage exploration. Building on the martingale approach of Jia and Zhou \cite{JiaZhou22a}, they derived a semi-analytical optimal Bernoulli distribution and established a policy improvement theorem. Related work on exploratory stopping under a singular control formulation similarly relies on cumulative residual entropy to regularize the randomized stopping time (see Dianetti et al. \cite{Dianettietal26}). These contributions share two structural features. First, they replace the stopping decision by a randomized policy (i.e., a stopping probability or a Bernoulli distribution over “stop” and “continue”) and analyze the resulting exploratory problem. Second, they introduce an entropy regularization term whose temperature controls the trade-off between exploration and bias; the optimal randomized policy is characterized through DPP. The entropy term is essential to their design because it keeps the randomized policy interior and yields a well-posed control problem, but it also introduces a temperature parameter whose tuning affects both convergence and the bias relative to the original stopping problem.
	
	The current paper develops an alternative model-free framework for continuous-time optimal stopping that does not rely on stochastic policy and entropy regularization. Instead, we work directly with the stochastic maximum principle (SMP) as in \cite{BoHuangWang26} for the data-driven stochastic control problem. The central idea is to replace the stopping problem by an equivalent intensity control problem under its randomized stopping framework of Gy\"ongy and $\check{\text{S}}$i$\check{\text{s}}$ka \cite{GyongySisika}. Rather than choosing a stopping probability at each instant, the agent controls the intensity (hazard rate) of a default-like stopping clock. This intensity formulation is value-function equivalent to the original optimal stopping problem in the sense that the optimal intensity encodes the same continuation/stopping regions as the classical variational inequality, but it lives naturally in the context of stochastic control and admits a maximum-principle treatment.
	
	The contributions of this paper are fourfold. First, we derive a performance-difference formula for the intensity-controlled problem (Proposition~\ref{prop:performancediff}), which compares any two admissible intensities through the stopping gap $\Delta^{\lambda}(t,x)=g(t,x)-v^{\lambda}(t,x)$, i.e., the discrepancy between the terminal payoff and the value under the current intensity. This yields a policy-improvement principle: if the candidate intensity is chosen so that the sign condition on the stopping gap holds, the candidate policy weakly improves the current one. Based on this principle, we propose an exact policy iteration for the intensity and prove its pointwise convergence to an optimal intensity $\lambda^*$ with $\mathcal J(\lambda^*)=\mathcal V^M$, together with an explicit factorial convergence rate (Theorem~\ref{thm:policy_convergence}). This provides an improvement mechanism that does not invoke entropy regularization. 
	
	Second, using the same performance-difference identity, we show that the stopping gap function $\Delta^{\lambda}(t,x)$ corresponding to an intensity strategy $\lambda$ can be learned from data by minimizing a weighted $L^2$-loss over a parameterized continuous function $h^{\theta}(t,x)$. Although $\Delta^{\lambda}$ also appears in the gradient of the weighted loss function, the first-order variation formula in Lemma~\ref{lem:DirectDerivJ} lets us express the loss gradient through the directional derivative $D\mathcal{J}(\lambda)[\nabla_{\theta}h^{\theta}]$, so that it can be optimized from realized trajectories and the candidate intensity alone (Lemma~\ref{lem:gradientStheta}). This way of expressing the loss gradient through a directional derivative is inspired by the idea of learning the Hamiltonian gradient in \cite{BoHuangWang26}. Its minimizer delivers an estimate of the stopping gap sufficient for policy improvement. This is where the model-free character of the framework becomes operational: neither the drift nor the volatility of the diffusion needs to be known. 
	
	Third, we use the learned stopping gap to drive a policy iteration on the $M$-capped intensity and establish its theoretical convergence. Because the associated Hamiltonian is linear in the intensity, the update is of bang-bang type and no entropy regularization is invoked. Combined with the learned gap, Theorem~\ref{prop:relaxed-policy-convergence} gives an explicit learning-error decomposition $\mathcal V^M-\mathcal J(\lambda^n)\le\mathcal E_M^n$, showing that the sub-optimality vanishes provided that the learning SMP residual $\mathcal R_n$ and the loss $S_M^{\lambda^n}(\theta^*)$ vanish along the iteration, and hence $\mathcal J(\lambda^n)\to\mathcal V^M$ as $n\to\infty$. 
	
	Fourth, we convert the learned intensity back into a genuine (approximate) optimal stopping time with respect to the original Brownian filtration of the prime stopping problem. To the best of our knowledge, this is the first RL algorithm in the continuous-time optimal stopping literature that passes from a randomized stopping rule/intensity to an (approximate) original-filtration stopping time for the original stopping problem. Given the terminal intensity $\lambda^{N^*}(t,x)$ via the iterative learning based on policy improvement result, we take a maximizer of the induced mapping $y\to\mathcal Q^{\lambda^{N^*}}(y)$ and define the original-filtration stopping time as in \eqref{eq:optimal-nearBST}. The error decomposition \eqref{eq:error_decomp} then shows that the resulting stopping time is a Brownian near-optimal stopping time with an explicit error bound. For implementation, the required gradient is given in Lemma~\ref{lem:learnQy}, the unobservable $\mu(t,x)$ is learned via the martingale orthogonality characterization proposed in Jia and Zhou \cite{JiaZhou22a}, and the maximizer is obtained by (inexact) gradient ascent. This yields an implementable RL procedure returning a genuine original-filtration stopping time rather than a randomized stopping rule.
	
	The remainder of the paper is organized as follows. Section \ref{sec:problem} formulates the model-free optimal stopping problem in continuous time. Section \ref{sec:randomizedformulation} introduces the randomized stopping based on stopping intensity control and the SMP is used to explore the capped intensity control problem. Section \ref{sec:policyimprovement} establishes a policy improvement theorem based on SMP and then learns the unknown stopping gap under a weighted $L^2$-loss. Section \ref{sec:learnBrownianstop} uses the trained stopping gap and the policy improvement theorem to learn the Brownian near-optimal stopping time of the prime problem. Numerical experiments based on the proposed algorithm are provided in Section \ref{sec:numerical}. The proofs of all results are provided in Section \ref{sec:proof_main}.
	
	{\bf Notations.}\quad We list below some notations that will be used frequently throughout the paper:
	\vspace{-0.6cm}
	\begin{center}
		\begin{longtable}{l l}
			$|\cdot|$ & Euclidean norm on Euclidean spaces\\
			$\nabla_x$ & Gradient operator with respect to $x\in\R^d$\\
			${\rm tr}(A)$ & Trace of a square matrix $A$, i.e., the sum of its diagonal elements\\
			$\E^{\Pb}$ & Expectation under $\Pb$ and we abbreviate it by $\E$ when $\Pb$ is clear from the context\\
			$\|X\|_{L^p}$ & $\|X\|_{L^p}:=(\E[\int_0^T|X_s|^p\d s])^{1/p}$ for a $\R^d$-valued process $X=(X_t)_{t\in[0,T]}$
		\end{longtable}
	\end{center} 
	
	\section{Problem Formulation}\label{sec:problem}
	
	In this section, we formulate the continuous time reinforcement learning (RL) problem for optimal stopping in a model-free setting. The objective of the agent is to learn an optimal (approxiamte) Brownian stopping time of continuous-time optimal stopping problem with unknown coefficients via stochastic maximum principle (SMP).
	
	Let $T>0$ be a finite terminal horizon and $d,r\in\mathbb{N}$. Consider a complete probability space $(\Omega^1,\F^1,\Pb^1)$ supporting a standard $r$-dimensional Brownian motion (BM) $W=(W_t)_{t\in[0,T]}$ and a $\R^d$-valued r.v. $\eta$ with $\E[|\eta|^2]<\infty$ which is independent of $W$. Furthermore, let $(\Omega^2,\F^2,\Pb^2)$ be another complete probability space supporting an exponentially distributed r.v. $\xi\sim\mathrm{Exp}(1)$. We construct the following product probability space $(\Omega,\F,\Pb)$ by setting $\Omega:=\Omega^1\times\Omega^2$, $\F:=\F^1\otimes\F^2$ and $\Pb:=\Pb^1\times\Pb^2$. For notational simplicity, we identify $\eta$, $W$ and $\xi$ with their respective extensions to $(\Omega,\F,\Pb)$. We then introduce the filtrations $\Fb^1=(\F_t^1)_{t\in[0,T]}$ and $\Fb=(\F_t)_{t\in[0,T]}$ on $\Omega$ respectively by $\F_t^1=\sigma(\eta,W_s;\,s\leq t)$ and $\F_t=\sigma(\eta,\xi,W_s;\,s\leq t)$ for all $t\in[0,T]$. Clearly, $\Fb^1$ is a subfiltration of $\Fb$ and $\xi$ is $\Pb$-independent of $(\eta,W)$. Without loss of generality, we may assume that both filtrations $\Fb^1$ and $\Fb$ satisfy the usual conditions. 
	
	Let the $\R^d$-valued state process $X=(X_t)_{t\in [0,T]}$ evolve as the following stochastic differential equation (SDE), for $t\in[0,T]$,
	\begin{align}\label{eq:SDE_dynamic}
		\d X_t=b(t,X_t)\d t+\sigma(t,X_t)\d W_t,\quad X_0=\eta,
	\end{align}
	where, $(b,\sigma):[0,T]\times\R^d\to\R^d\times\R^{d\times r}$ are Borel measurable functions so that Eq.~\eqref{eq:SDE_dynamic} admits a unique (strong, i.e., $\Fb^1$-adapted) solution (see Assumption \ref{ass1}-(A1) below). Let ${\cal T}$ be the set of $\Fb^1$-stopping times taking values in $[0,T]$. The goal is to maximize the following expected reward process over all $\tau\in{\cal T}$:
	\begin{align}\label{eq:reward_func}
		J(\tau):=\E\left[\Phi_{\tau}\right]=\E\left[\int_0^{\tau} f(s,X_s)\d s+g(\tau,X_{\tau})\right]\longrightarrow V:=\sup_{\tau\in{\cal T}}J(\tau),
	\end{align}
	where, the running reward $f:[0,T]\times\R^d\to\R$ and the terminal reward $g:[0,T]\times\R^d\to\R$ are measurable mappings.
	
	In the model-free RL framework, the agent has {\it no access} to the model coefficients $(b,\sigma,f,g)$. Instead, the agent interacts with the environment \eqref{eq:SDE_dynamic} by exploring different policies and improves her policy through the observed states and rewards. More precisely, for an arbitrary time discretization given by $0=t_0<t_1<\cdots<t_{N_T}=T$, the agent can only observe $N$ independent realizations of the state-reward-stopping given by $\{X_{t_j},f(t_j,X_{t_j}),g(t_j,X_{t_j})\}_{j=0}^{N_T}$. Hence, the cumulative reward at observable time $t_j$ is given by $\Phi_{t_j}=\sum_{i=1}^jf(t_i,X_{t_i})(t_i-t_{i-1})+g(t_j,X_{t_j})$ for $j=0,1,\ldots,N_T$. Correspondingly, we denote these samples by $\{(X_{t_j,i}, f_{t_j,i},g_{t_j,i})\}_{i=1}^{N}$ for $j=0,1,\ldots,N_T$.
	
	Lastly, we impose the following assumptions on model coefficients throughout the paper.
	\begin{ass}\label{ass1}
		\begin{itemize}
			\item [{\rm (A1)}] The functions $(b,\sigma):[0,T]\times\R^d\to \R^d\times\R^{d\times r}$ are Lipschitzian continuous with linear growth with respect to $x\in\R^d$ uniformly in $t\in[0,T]$.
			
			\item [{\rm (A2)}] The running function $f:[0,T]\times\R^d\to\R$ is continuous and has the squared growth w.r.t. $x\in\R^d$ uniformly in $t\in[0,T]$, and the terminal function $g:[0,T]\times\R^d\to\R$ is Lipschitzian continuous in $(t,x)\in[0,T]\times\R^d$. 
		\end{itemize}
	\end{ass}
	
	\section{Randomized Stopping Formulation}\label{sec:randomizedformulation}   
	We now consider the optimal stopping problem \eqref{eq:reward_func} ‌in its randomized form‌ under the framework of Gy\"ongy and $\check{\text{S}}$i$\check{\text{s}}$ka \citep{GyongySisika}. To this end, let $\mathscr{U}[0,T]$ be the collection of all nonnegative  $\Fb^1$-progressively measurable processes $\lambda=(\lambda_t)_{t\in [0,T]}$. For any $\lambda=(\lambda_t)_{t\in[0,T]}\in\mathscr{U}[0,T]$, we introduce the $\Fb$-stopping time that
	\begin{align}\label{eq:taulambda}
		\tau_{\xi}^{\lambda}:=\inf\left\{t\geq 0;~\int_0^t\lambda_s\d s\geq\xi\right\}\wedge T,   
	\end{align}
	where, the convention $\inf\varnothing=+\infty$. 
	
	Note that $\tau_{\xi}^{\lambda}$ is not an $\Fb^1$-stopping time, and hence $\tau_{\xi}^{\lambda}\notin{\cal T}$. To establish an $\Fb^1$-stopping time related to $\tau_{\xi}^{\lambda}$, a natural candidate is to introduce a counterpart of $\tau_{\xi}^{\lambda}$ at an arbitrary sampling point of the random (standard) exponential threshold $\xi$. In other words, we define the following parameterized $\Fb^1$-stopping time by, for any $y>0$, 
	\begin{align}\label{eq:tauylambda}
		\tau_y^{\lambda}:=\inf\left\{t\geq 0;~\int_0^t\lambda_s\d s\geq y\right\}\wedge T.
	\end{align}
	Obviously, $\tau_y^{\lambda}\in{\cal T}$ for any $y>0$. Furthermore, let us define that, for any $y>0$,
	\begin{align}\label{eq:calQy}
		{\cal Q}^{\lambda}(y):=J(\tau_{y}^{\lambda})=\E\left[\Phi_{\tau_y^{\lambda}}\right],  
	\end{align}
	where, the cumulative reward process $\Phi=(\Phi_t)_{t\in[0,T]}$ is defined in \eqref{eq:reward_func}. Then, we have the following result.
	\begin{lemma}\label{lem:valuefcnrep}
		For any $\lambda\in\mathscr{U}[0,T]$, it holds that 
		\begin{align}\label{eq:vaefcntrp0}
			\mathcal{J}(\lambda)&:=\E\left[\Phi_{\tau_{\xi}^{\lambda}}\right]=\int_0^{\infty}{\cal Q}^{\lambda}(y)e^{-y}\d y\nonumber\\
			&=\E\left[\int_0^TL_s^{\lambda}(f(s,X_s)+\lambda_sg(s,X_s))\d s+L_T^{\lambda}g(T,X_T)\right].    
		\end{align}
		Here, $L^{\lambda}=(L_t^{\lambda})_{t\in[0,T]}$ is the survival process defined by
		\begin{align}\label{eq:exponenclock}
			L_t^{\lambda}:=\exp\left(-\int_0^t\lambda_r\d r\right),\quad\forall t\in[0,T]. 
		\end{align}
	\end{lemma}
	
	Lemma \ref{lem:valuefcnrep} motivates us to consider the following optimal stopping intensity problem described as:
	\begin{align}\label{eq:optimalintensity}
		{\cal V}:=\sup_{\lambda\in\mathscr{U}[0,T]}\mathcal{J}(\lambda).
	\end{align}
	By virtue of Theorem 3.2 in Gy\"ongy and $\check{\text{S}}$i$\check{\text{s}}$ka \cite{GyongySisika}, we have the equivalence between the original optimal stopping problem \eqref{eq:reward_func} and the optimal intensity problem \eqref{eq:optimalintensity} in the sense that $V={\cal V}$. We next consider the optimal control problem \eqref{eq:optimalintensity} subject to the underlying state process given by, for $t\in[0,T]$,
	\begin{align}\label{eq:jointstateXL}
		\begin{cases}
			\displaystyle \d X_t = b(t,X_t)\d t + \sigma(t, X_t)\d W_t,~X_0=\eta;\\[0.4em]
			\displaystyle \d L_t^{\lambda} = -L_t^{\lambda}\lambda_t\d t,~L_0^{\lambda}=1.
		\end{cases}    
	\end{align}
	The Hamiltonian $H:[0,T]\times\R^d\times(0,\infty)\times[0,\infty)\times\R^{d+1}\times\R^{(d+1)\times r}\to\R$ corresponding to problem \eqref{eq:optimalintensity}-\eqref{eq:jointstateXL} is given by
	\begin{align}\label{eq:HM}
		H\left(t,x,\ell,a,\begin{pmatrix}  y\\ p\end{pmatrix},\begin{pmatrix} z\\ q\end{pmatrix}\right)&:= b(t,x)^{\top}y-\ell a p + {\rm tr}\left(\begin{pmatrix} \sigma(t,x)\\ 0_{1\times r}\end{pmatrix}^{\top}\begin{pmatrix} z\\ q\end{pmatrix}\right)\nonumber\\
		&\quad + \ell(f(t,x)+a g(t,x)).
	\end{align}
	For any admissible intensity control process $\lambda\in\mathscr{U}[0,T]$, consider the following adjoint BSDE given by, for $t\in[0,T]$,
	\begin{align}\label{eq:BSDE}
		\begin{cases}
			\displaystyle  \d\begin{pmatrix}  Y_t^{\lambda}\\ P_t^{\lambda}\end{pmatrix}=-\begin{pmatrix}   \nabla_{x} H\left(t,\begin{pmatrix} X_t\\ L_t^{\lambda}\end{pmatrix},\lambda_t,\begin{pmatrix}Y_t^{\lambda}\\ P_t^{\lambda}\end{pmatrix},\begin{pmatrix}Z_t^{\lambda}\\ Q_t^{\lambda}\end{pmatrix}\right) \\ f(t,X_t)+\lambda_t(g(t,X_t)-P_t^{\lambda})\end{pmatrix}\d t+\begin{pmatrix}  Z_t^{\lambda}\\ Q_t^{\lambda}\end{pmatrix}\d W_t,\\[2em]
			\displaystyle ~~\begin{pmatrix}  Y_T^{\lambda}\\ P_T^{\lambda}\end{pmatrix}=\begin{pmatrix} L_T^{\lambda}\nabla_{x}g(T,X_T)\\ g(T,X_T)\end{pmatrix}.
		\end{cases}
	\end{align}
	
	Let $\mathscr{U}^{\rm m}[0,T]$ be the set of all admissible stopping intensity $\lambda\in\mathscr{U}[0,T]$ with the Markovian (feedback) form $\lambda_t=\lambda(t,X_t)$ for $t\in[0,T]$ for some measurable feedback function $(t,x)\to\lambda(t,x)$ on $[0,T]\times\R^d$ satisfying the polynomial growth in $x\in\R^d$ uniformly in $t\in[0,T]$. Then, we have the following decoupling field of the adjoint process of the survival process.
	\begin{lemma}[Decoupling Field of BSDE]\label{lem:BSDElambdadecoupling}
		For any $\lambda\in\mathscr{U}^{\rm m}[t,T]$, there exists a continuous decoupling field $v^{\lambda}:[0,T]\times\R^d\to\R$ such that 
		\begin{align}\label{eq:decouplfields}
			P_t^{\lambda}=v^{\lambda}(t,X_t),\quad \forall t\in[0,T].
		\end{align}
		Furthermore, this decoupling field admits the following probabilistic representation given by, for $(t,x)\in[0,T]\times\R^d$,
		\begin{align}\label{eq:probrep}
			v^{\lambda}(t,x)&=\E_{t,x}\bigg[\int_t^Te^{-\int_t^r\lambda(s,X_s)\d s}\left(f(r,X_r)+\lambda(r,X_r)g(r,X_r)\right)\d r\nonumber\\
			&\qquad\qquad+e^{-\int_t^T\lambda(s,X_s)\d s}g(T,X_T)\bigg]. 
		\end{align}
		Here, the conditional expectation operator $\E_{t,x}[\cdot]:=\E[\cdot|X_t=x]$. Consequently, we have
		\begin{align}\label{eq:vlambdaJ}
			{\cal J}(\lambda)=\E\left[v^{\lambda}(0,X_0)\right].   
		\end{align}
	\end{lemma}
	
	We observe that the Hamiltonian $H(\cdot)$ is linear in the action variable $a$ (see Eq.~\eqref{eq:HM}). This implies that the optimal stopping intensity would be a bang-bang control taking values $0$ and $\infty$. Hence, it is not appropriate to use this control to learn the original optimal stopping time. To overcome this difficulty, we introduce a truncated version of problem \eqref{eq:optimalintensity} given by, for any $M>1$, 
	\begin{align}\label{eq:optimalintensityM}
		{\cal V}^M:=\sup_{\lambda\in\mathscr{U}^M[0,T]}\mathcal{J}(\lambda),
	\end{align}
	where, $\mathscr{U}^M[0,T]$ is the set of admissible stopping intensities in $\mathscr{U}[0,T]$ taking values in $[M^{-1},M]$. By standard moment estimates, we can verify that $\mathcal{V}^M<\infty$.
	
	Obviously, $M\to{\cal V}^M$ is increasing and ${\cal V}^M\leq{\cal V}$ for all $M>0$. The following lemma establishes the convergence of the optimal value ${\cal V}^M$ to the optimal value ${\cal V}$ as $M$ tends to infinity.
	\begin{lemma}[Convergence of Truncated Optimal Values]\label{lem:VM-V-convergence}
		For any $M>1$, there exists a constant $C>0$ independent of $M$ such that
		\begin{align}\label{eq:estimate2}
			&0\leq {\cal V}-{\cal V}^M\leq C\left(\frac{1}{M^2}+\frac{1}{M} + \frac{1}{\sqrt{M}}\right).
		\end{align}
	\end{lemma}
	The following lemma identifies an optimal stopping intensity in $\mathscr{U}^M[0,T]$ for the truncated problem \eqref{eq:optimalintensityM}.
	\begin{lemma}[Optimal Truncated Stopping Intensity]\label{lem:cappedintensity} Let $M>1$. Consider the following scalar BSDE satisfied by $(P^M,Q^M)=(P_t^M,Q_t^M)_{t\in[0,T]}$ which is given by, for $t\in[0,T]$,
		\begin{align}\label{eq:BSDEPmM}
			P_t^{M}&=g(T,X_T)+\int_t^T\left\{f(s,X_s)+M\left(g(s,X_s)-P_s^{M}\right)^+ - \frac{1}{M}\left(g(s,X_s)-P_s^{M}\right)^{-}\right\}\d s\nonumber\\
			&\quad-\int_t^T Q_s^{M}\d W_s.    
		\end{align}
		Then, for any $\lambda\in\mathscr{U}^{M}[0,T]$, we have $P_t^{\lambda}\leq P_t^{M}$ for all $t\in[0,T]$, a.s.. Furthermore, let us define that
		\begin{align}\label{eq:lambdaMstar}
			\lambda_t^{M,*}:=M\mathbf{1}_{\{g(t,X_t)\geq P_t^M\}}+\frac{1}{M}\mathbf{1}_{\{g(t,X_t)< P_t^M\}},\quad \forall t\in[0,T].
		\end{align}
		Then, $\lambda^{M,*}=(\lambda_t^{M,*})_{t\in[0,T]}\in\mathscr{U}^M[0,T]$ is an optimal stopping intensity  of the truncated problem $\sup_{\lambda\in\mathscr{U}^M[0,T]}{\cal J}(\lambda)$.
	\end{lemma}
	
	Let $\mathscr{U}^{M,{\rm m}}[0,T]:=\mathscr{U}^M[0,T]\cap\mathscr{U}^{\rm m}[0,T]$. Then, we introduce the following dynamical version of problem \eqref{eq:optimalintensity}, for $(t,x)\in [0,T]\times\R^d$,
	\begin{align}\label{eq:optimalintensityMfeedback}
		&{\cal V}^M(t,x):=\sup_{\lambda\in\mathscr{U}^{M,{\rm m}}[t,T]}v^{\lambda}(t,x)=v^{\lambda^{M,*}}(t,x)
	\end{align}
	with the decoupling field function $v^{\lambda}(t,x)$ being given by \eqref{eq:probrep}. It follows from Lemma \ref{lem:BSDElambdadecoupling} and Lemma \ref{lem:cappedintensity} that ${\cal V}^M(t,x)$ is a decoupling field of the solution $P^{M}=(P_t^M)_{t\in[0,T]}$ of the BSDE \eqref{eq:BSDEPmM}, i.e., $P_t^M={\cal V}^M(t,X_t)$ for all $t\in[0,T]$. Hence, the optimal stopping intensity \eqref{eq:lambdaMstar} admits the following representation, for any $(t,x)\in[0,T]\times\R^d$,
	\begin{align}\label{eq:lambdaMstardecoupling}
		\lambda^{M,*}(t,x):=M\mathbf{1}_{\{g(t,x)-{\cal V}^M(t,x)\geq 0\}}+\frac{1}{M}\mathbf{1}_{\{g(t,x)-{\cal V}^M(t,x)<0\}}.
	\end{align}
	The optimal truncated stopping intensity above motivates the establishment of the policy improvement result, which is discussed in the next section.
	
	\section{Policy Improvement}\label{sec:policyimprovement}
	
	To establish the policy improvement based on SMP, we first have the following result on the performance difference of the intensity control problem:
	\begin{prop}[Performance Difference]\label{prop:performancediff}
		For any $\lambda\in\mathscr{U}^{\rm m}[0,T]$ and any $\Fb^1$-adapted real-valued process $\phi=(\phi_t)_{t\in [0,T]}$ with a lower bound satisfying $\|\phi\|_{L^2}<\infty$, let $(P^{\lambda},Q^{\lambda})=(P_t^{\lambda},Q_t^{\lambda})_{t\in[0,T]}$ be the adjoint process of the survival process satisfying BSDE \eqref{eq:BSDE}. Then, it holds that
		\begin{align}\label{eq:differenceJ02}
			P_t^{\phi}-P_t^{\lambda}&= \E\left[\int_t^T(L_t^{\phi})^{-1}L_s^{\phi}(\phi_s-\lambda_s)\Delta^{\lambda}(s,X_s)\d s\Big|\F_t^1\right],\quad \forall t\in [0,T].  
		\end{align}
		Here, the survival process $L_s^{\phi}=\exp(-\int_0^s\phi_r\d r)$ for $s\in[0,T]$, and the {\it stopping gap} $\Delta^{\lambda}(t,x)$ is defined by, for any $(t,x)\in[0,T]\times\R^d$,
		\begin{align}\label{eq:stopgap}
			\Delta^{\lambda}(t,x):=g(t,x)-v^{\lambda}(t,x),  
		\end{align}
		where, the value function $v^{\lambda}(t,x)$ for $(t,x)\in[0,T]\times\R^d$ is given by \eqref{eq:probrep}. In particular, we have
		\begin{align}\label{eq:stopgap2}
			\mathcal{J}(\phi)-\mathcal{J}(\lambda)&= \E\left[\int_0^TL_s^{\phi}(\phi_s-\lambda_s)\Delta^{\lambda}(s,X_s)\d s\right]. 
		\end{align}
	\end{prop}
	
	The following lemma establishes the stability of stopping gap process $\Delta_t^\lambda:=\Delta^{\lambda}(t,X_t)$ for $t\in[0,T]$ under the intensity process $\lambda=(\lambda_t)_{t\in[0,T]}\in\U^{M,\rm m}[0,T]$. Its proof can follow from standard moment estimates and the stability theory of BSDEs, and is hence omitted.
	\begin{lemma}[Stability of BSDE]\label{lem:delta_stability}
		Let
		$\lambda^1,\lambda^2\in\U^{M,\rm m}[0,T]$, $(P^{\lambda^1},Q^{\lambda^1})$ and $(P^{\lambda^2},Q^{\lambda^2})$ be the corresponding adjoint processes defined in BSDE \eqref{eq:BSDE}. Then, there exists a constant $C>0$ depending on $(M,T)$ only such that
		\begin{align*}
			\mathbb E\left[
			\sup_{t\in [0,T]}\left|\Delta_t^{\lambda^1}-\Delta_t^{\lambda^2}\right|^2\right]=\mathbb E\left[
			\sup_{t\in [0,T]}\left|P_t^{\lambda^1}-P_t^{\lambda^2}\right|^2\right]
			\leq C\left\|\lambda^1-\lambda^2\right\|_{L^2}^2.
		\end{align*}
	\end{lemma}
	
	The following lemma also explains that the stopping gap $\Delta^{\lambda}(t,x)$ defined by \eqref{eq:stopgap} is related to the optimal stopping of the prime problem.
	\begin{lemma}\label{lem:stoptimedelta}
		For any $\lambda\in\mathscr{U}^{{\rm m}}[0,T]$ and the stopping gap $\Delta^{\lambda}(t,x)$ given by \eqref{eq:stopgap}, let us define the following $\Fb^1$-stopping time by  
		\begin{align}\label{eq:stoptimetaulambda}
			\tau_{\lambda}:= \inf\left\{t\geq0;~\Delta^{\lambda}(t,X_t)\geq0\right\}\wedge T.   
		\end{align}
		Then, $V\geq J(\tau_{\lambda})\geq{\cal J}(\lambda)$. Furthermore, set $\delta_{\lambda}:={\cal V}-{\cal J}(\lambda)$. Hence, $\tau_{\lambda}\in{\cal T}$ is a $\delta_{\lambda}$-optimal stopping time of the original problem \eqref{eq:reward_func}. In particular, if $\lambda^*\in\mathscr{U}^{{\rm m}}[0,T]$ is the optimal stopping intensity in the sense that ${\cal V}={\cal J}(\lambda^*)$, then $\tau_{\lambda^*}\in{\cal T}$ is an optimal stopping time of the original problem \eqref{eq:reward_func}.
	\end{lemma}
	
	In particular, set $\Delta^M(t,x):=g(t,x)-v^{\lambda^{M,*}}(t,x)=g(t,x)-{\cal V}^M(t,x)$ for $(t,x)\in[0,T]\times\R^d$. Here, the bounded stopping intensity is given by $\lambda_t^{M,*}=\lambda^{M,*}(t,X_t)=M{\bf1}_{\{\Delta^M(t,X_t)\geq0\}}+M^{-1}{\bf1}_{\{\Delta^M(t,X_t)<0\}}$ for $t\in[0,T]$ (see also \eqref{eq:lambdaMstardecoupling}). Similarly to Lemma \ref{lem:stoptimedelta}, the capped stopping gap $\Delta^M(t,x)$ is also related to the near-optimal stopping of the prime problem.
	\begin{corollary}[Near-Optimal Brownian Stopping Time]\label{coro:approximatestoptime}
		For any $M>1$, let $\tau_{\lambda^{M,*}}:=\inf\{t\geq0;~\Delta^{M}(t,X_t)\geq0\}$. Then, $\tau_{\lambda^{M,*}}\in{\cal T}$ is a $\varepsilon_M$-optimal stopping time of the prime problem \eqref{eq:reward_func}. Here,  $\varepsilon_M=C(M^{-2}+{M}^{-1}+M^{-\frac{1}{2}})$ is given in Lemma \ref{lem:VM-V-convergence}. 
	\end{corollary}
	
	Based on Proposition~\ref{prop:performancediff}, we can establish the following policy improvement result.
	\begin{lemma}[Policy Improvement]\label{lem:policy_improvement}
		Let $M>1$. For any $\lambda\in\mathscr{U}^{M,{\rm m}}[0,T]$, define the intensity process $\hat{\lambda}^{M}=(\hat{\lambda}_t^{M})_{t\in [0,T]}\in\mathscr{U}^{M,{\rm m}}[0,T]$ by 
		\begin{align*}
			\hat{\lambda}_t^{M}:=M\mathbf{1}_{\{\Delta^{\lambda}(t,X_t)\geq0\}} + \frac{1}{M} \mathbf{1}_{\{\Delta^{\lambda}(t,X_t)<0\}},\quad \forall t\in[0,T].  
		\end{align*}
		Then, we have the policy improvement that $\mathcal{J}(\hat{\lambda}^M)\geq \mathcal{J}(\lambda)$.
	\end{lemma}
	
	We now state the first main result of this section: the policy iteration based on Lemma~\ref{lem:policy_improvement} is globally convergent; that is, for any initial policy $\lambda\in\U^{M,\rm m}[0,T]$, the policy iteration converges to an optimal intensity policy in $\U^{M,\rm m}[0,T]$ after sufficiently many iterations.
	\begin{theorem}[Policy Iteration Convergence]\label{thm:policy_convergence}
		For any initial intensity $\lambda^0=\lambda\in\mathscr{U}^{M,m}[0,T]$, let the $n$-th policy iteration of intensity be given by, for $n\geq0$,
		\begin{align}\label{eq:iterative000}
			\lambda_t^{n+1}= M\mathbf{1}_{\{\Delta^{\lambda^n}(t,X_t)\geq0\}} + \frac{1}{M} \mathbf{1}_{\{\Delta^{\lambda^n}(t,X_t)<0\}},\quad \forall t\in[0,T],
		\end{align}
		where, $\Delta^{\lambda^n}(t,x)=g(t,x)-v^{\lambda^n}(t,x)$ is the stopping gap defined by \eqref{eq:stopgap} w.r.t. the $n$-th iterative intensity. Then, as $n\to\infty$, the intensity $\lambda_t^n:=\lambda^n(t,X_t)$ converges to some optimal intensity $\lambda^*_t=\lambda^*(t,X_t)$, $\mathbb{P}$-a.s. in the sense that $\mathcal{J}(\lambda^*)=\mathcal{V}^M$. Moreover, there exists a constant $C>0$ independent of $n$ such that, for any $n\geq1$,
		\begin{align}\label{eq:converge_rate}
			\left|\mathcal{J}(\lambda^n)-\mathcal{V}^M\right|\leq C\frac{(MT)^{n}}{n!}.
		\end{align} 
	\end{theorem}
	
	\begin{remark}
		Let $\lambda^*=(\lambda_t^*)_{t\in[0,T]}$ be the limit point defined in Theorem~\ref{thm:policy_convergence}. Then, by virtue of Corollary~\ref{coro:approximatestoptime}, the induced $\Fb^1$-stopping time $\tau_{\lambda^*}=\inf\{t\geq 0;\Delta^{\lambda^*}(t,X_t)\geq 0\}\wedge T$ is an $\varepsilon_M$-optimal solution to the prime problem \eqref{eq:reward_func}.
	\end{remark}
	\begin{remark}
		Theorem \ref{thm:policy_convergence} above provides a theoretical motivation for our RL algorithm (Algorithm~\ref{alg:optimal_stopping} in Section \ref{sec:learnBrownianstop}). In particular, it shows that iteratively improving the policy according to the stopping gap function leads to the convergence to an optimal intensity policy, with a factorial convergence rate. This suggests that  one can focus on learning the stopping gap function and then iteratively update the intensity policy based on the learned stopping gap.  We further establish the corresponding convergence of the resulting algorithm in Theorem~\ref{prop:relaxed-policy-convergence}.
	\end{remark}
	
	Proposition \ref{prop:performancediff} can also help us to establish the directional derivative of the objective function ${\cal J}(\lambda)$ directly, which is given by the following lemma whose proof is omitted:
	\begin{lemma}\label{lem:DirectDerivJ}
		For any stopping intensity $\lambda\in\mathscr{U}^{{\rm m}}[0,T]$ and any $\Fb$-adapted real-valued process $\beta=(\beta_t)_{t\in [0,T]}$ satisfying $\|\beta\|_{L^2}<\infty$, we have the first-order variation equality given by
		\begin{align}\label{eq:DJbeta}
			\frac{\mathcal{J}(\lambda+\epsilon\beta)-\mathcal{J}(\lambda)}{\epsilon}=\E\left[\int_0^T L_s^{\lambda}\beta_s\Delta^{\lambda}(s,X_s)\d s\right]+O(\epsilon\|\beta\|_{L^2}).
		\end{align}
	\end{lemma}
	
	
	Let $\lambda\in\mathscr{U}^{M,{\rm m}}[0,T]$ and the stopping gap $\Delta^{\lambda}(t,x)$ be given by \eqref{eq:stopgap}. By virtue of the policy improvement result (Lemma \ref{lem:policy_improvement}), we next learn the stopping gap $\Delta^{\lambda}(t,x)$ by using the parameterized continuous function $h^{\theta}:[0,T]\times\R^d\to\R$ with the unknown parameter $\theta\in\R^q$ such that $\theta\to h^{\theta}(t,x)$ is $C^1(\R^q)$. We then introduce the following weighted $L^2$-loss function given by
	\begin{align}\label{eq:learnDeltaloss}
		S_M^{\lambda}(\theta):=\frac{1}{2}\E\left[\int_0^T L_s^{\lambda}\left|h^{\theta}(s,X_s)-\Delta^{\lambda}(s,X_s)\right|^2\d s\right],   
	\end{align}
	where the weight process $L_t^{\lambda}:=\exp(-\int_0^t\lambda_s\d s)$ for $t\in[0,T]$. The gradient of the loss $S_M^{\lambda}(\theta)$ with respect to $\theta\in\R^q$ is given by
	\begin{align*}
		\nabla_{\theta}S_M^{\lambda}(\theta)&=\E\left[\int_0^T L_s^{\lambda}(h^{\theta}(s,X_s)-\Delta^{\lambda}(s,X_s))\nabla_{\theta}h^{\theta}(s,X_s)\d s\right]\nonumber\\
		&=\E\left[\int_0^T L_s^{\lambda}h^{\theta}(s,X_s)\nabla_{\theta}h^{\theta}(s,X_s)\d s\right]-\E\left[\int_0^T L_s^{\lambda}\Delta^{\lambda}(s,X_s)\nabla_{\theta}h^{\theta}(s,X_s)\d s\right].
	\end{align*}
	To use the gradient descent to learn the parameter $\theta\in\R^q$, we need to approximate the second term of the last line in the above display (since $\Delta^{\lambda}$ is the function we want to learn) by using the directional derivative $D{\cal J}(\lambda)$ provided in Lemma \ref{lem:DirectDerivJ}. In fact, we have from Lemma \ref{lem:DirectDerivJ} that
	\begin{lemma}\label{lem:gradientStheta}
		The gradient $\nabla_{\theta}S_M^{\lambda}(\theta)$ can be represented by, for any $\theta\in\R^q$,
		\begin{align}\label{eq:nablathetaL}
			\nabla_{\theta}S_M^{\lambda}(\theta)&=\E\left[\int_0^T L_s^{\lambda}h^{\theta}(s,X_s)\nabla_{\theta}h^{\theta}(s,X_s)\d s\right]-D\mathcal{J}(\lambda)[\nabla_{\theta}h^{\theta}], 
		\end{align}    
		where $D\mathcal{J}(\lambda)[\nabla_{\theta}h^{\theta}]=(D\mathcal{J}(\lambda)[\partial_{\theta_j}h^{\theta}];j=1,\ldots,q)^{\top}$ is a $q$-dimensional column vector whose $j$-th component is the directional derivative of ${\cal J}(\lambda)$ along the direction $\beta=(\partial_{\theta_j}h^{\theta}(t,X_t))_{t\in[0,T]}$. 
	\end{lemma}

	\section{Learning Brownian Near-Optimal Stopping}\label{sec:learnBrownianstop}
	
	In this section, we propose the algorithm for learning Brownian near-optimal stopping time of the prime problem \eqref{eq:reward_func}.
	
	Let $M>1$ be fixed. For $n\geq0$, let $h^{\theta^*(\lambda^n)}(t,x)$ be an approximation to the stopping gap $\Delta^{\lambda^n}(t,x)=g(t,x)-v^{\lambda^n}(t,x)$ corresponding to the $n$-th iterative stopping intensity $\lambda^n\in\mathscr{U}^{M,{\rm m}}[0,T]$ under the weighted  $L^2$-loss $S_M^{\lambda^n}(\theta^*)$ given by \eqref{eq:learnDeltaloss}. We then consider the following  policy iteration:
	\begin{align}\label{eq:relaxed-policy-iteration}
		\lambda_t^{n+1}
		=(1-\rho_n)\lambda_t^n
		+\rho_n \left(M\mathbf{1}_{\{h^{\theta^*(\lambda^n)}(t,X_t)\geq0\}}+M^{-1}\boldsymbol{1}_{\{h^{\theta^*(\lambda^n)}(t,X_t)< 0\}}\right),\quad \forall t\in[0,T],
	\end{align}
	where $(\rho_n)_{n\geq0}$ is a hyperparameter sequence. In particular,  $\rho_n=1$ corresponds to the hard deployment, i.e.,  $\lambda^{n+1}_t=M\mathbf{1}_{\{h^{\theta^*(\lambda^n)}(t,X_t)\geq0\}}+M^{-1}\boldsymbol{1}_{\{h^{\theta^*(\lambda^n)}(t,X_t)< 0\}}$ for $t\in[0,T]$. 
	To measure the violation of the Hamiltonian maximization condition, let us define the following SMP residual by, for $\lambda\in\mathscr{U}^{M,{\rm m}}[0,T]$,
	\begin{align}\label{eq:SMP-residual}
		&\mathcal{R}^M(\lambda)\\
		&\quad:=\mathbb E\left[\int_0^TL_t^\lambda
		\left(\left(M-\lambda(t,X_t)\right)(\Delta^\lambda(t,X_t))^+ + (\lambda(t,X_t)-M^{-1})(\Delta^\lambda(t,X_t))^-\right)\d t\right].\nonumber
	\end{align}
	Proposition~\ref{prop:performancediff} yields that ${\cal R}^M(\lambda)=\sup_{\mu\in\mathscr{U}^{M,{\rm m}}[0,T]}D{\cal J}(\lambda)[\mu-\lambda]\geq0$. Thus, $\mathcal{R}^M(\lambda)=0$ iff $\lambda(t,X_t)=M$ when $\Delta^\lambda(t,X_t)\geq0$; while $\lambda(t,X_t)=M^{-1}$ when $\Delta^\lambda(t,X_t)<0$, $\d t\times\d\mathbb{P}$-a.s.. Thus, by using Proposition~\ref{prop:performancediff} again and the SMP, we conclude that $\lambda\in\U^{M,{\rm m}}[0,T]$ is optimal iff $\mathcal{R}^M(\lambda)=0$. 
	
	We now turn to the second main result of this paper, which provides a theoretical convergence guarantee together with an error analysis for our RL algorithm (Algorithm~\ref{alg:optimal_stopping}).
	\begin{theorem}[Learning Error Decomposition]\label{prop:relaxed-policy-convergence}
		Let $M>1$ be fixed, $\mathcal R_n^M:=\mathcal{R}^M(\lambda^n)$ for $n\geq0$, and the sequence $(\rho_n)_{n\geq0}$ with $\rho_n\in(0,1]$ satisfies $\sum_{n=0}^\infty\rho_n=\infty$ and $\sum_{n=0}^\infty
		\rho_n\sqrt{S_M^{\lambda^n}(\theta^*)}<\infty$. Then, $\lim_{n\to\infty}{\cal R}^M_n=0$ and  $\lim_{n\to\infty}\mathcal{J}(\lambda^n)=\mathcal{V}^M$. Furthermore, the following error decomposition holds, for any $n\geq0$,
		\begin{align}\label{eq:error_decomp}
			0&\leq  \mathcal{V}^M-\mathcal{J}(\lambda^n)\leq{\cal E}_M^n:=e^{MT}\mathcal{R}_n^M.
		\end{align}
	\end{theorem}
	
	\begin{remark}
		If we take $\rho_n=1$ for all $n\geq0$ and assume that the stopping gap is learned exactly at each iteration, namely $S_M^{\lambda^n}(\theta^*)=0$ with $n\geq0$, the learning framework reduces to the exact policy iteration setting considered in Theorem~\ref{thm:policy_convergence}.
	\end{remark}
	
	\begin{remark}
		Theorem~\ref{prop:relaxed-policy-convergence} proposes a method for learning the optimal intensity that is different from the approaches in Dong~\cite{Dong24} and Dai et al.~\cite{Daisunxuzhou26}, together with the corresponding theoretical convergence guarantee. A key difference is that, unlike \cite{Dong24,Daisunxuzhou26}, we shall discuss later how to construct from an intensity $\lambda=(\lambda_t)_{t\in[0,T]}$ an $\Fb^1$-stopping time $\tau^\lambda$ that improves the value function, i.e., $\mathcal J(\lambda)\leq J(\tau^\lambda)$. This is fundamentally different from the constructions in \cite{Dong24} and \cite{Daisunxuzhou26}: therein, once an intensity is obtained, one can only obtain an $\Fb$-stopping time $\tau^{\lambda}_{\xi}$ (i.e., a randomized stopping time), rather than a stopping time with respect to the original Brownian filtration.
	\end{remark}
	
	Let the outer-loop terminate at the $N^*$-th iteration with the optimal learning parameter $\theta^*(\lambda^{N^*})$ and the corresponding learning intensity function be $\lambda^{N^*}(t,x)$. Then, it follows from \eqref{eq:error_decomp} that $0\leq  \mathcal{V}^M-\mathcal{J}(\lambda^{N^*})\leq {\cal E}_M^{N^*}$, where ${\cal E}_M^{N^*}$ is the error term on the RHS of \eqref{eq:error_decomp} with $n={N^*}$. Recall the mapping $y\to{\cal Q}^{\lambda}(y)$ defined by \eqref{eq:calQy} for $\lambda\in\U[0,T]$. Assume that there exists a constant $\hat{y}>0$ such that $\sup_{y>0}{\cal Q}^{\lambda^{N^*}}(y)={\cal Q}^{\lambda^{N^*}}(\hat{y})$, i.e., $\hat{y}$ is a maximizer of $y\to{\cal Q}^{\lambda^{N^*}}(y)$ on $(0,\infty)$. Let us introduce the Brownian stopping time $\tau_{\hat{y}}{^{\lambda^{N^*}}}$ in terms of \eqref{eq:tauylambda}, i.e., 
	\begin{align}\label{eq:optimal-nearBST}
		\tau_{\hat{y}}{^{\lambda^{N^*}}}=\inf\left\{t\geq 0;~\int_0^t\lambda_s^{N^*}\d s\geq \hat{y}\right\}\wedge T.
	\end{align}
	Then, we have from Lemma \ref{lem:valuefcnrep} that
	\begin{align*}
		\mathcal{J}(\lambda^{N^*})=\E\left[\Phi_{\tau_{\xi}^{\lambda^{N^*}}}\right]=\int_0^{\infty}{\cal Q}^{\lambda^{N^*}}(y)e^{-y}\d y\leq \sup_{y>0}{\cal Q}^{\lambda^{N^*}}(y)={\cal Q}^{\lambda^{N^*}}(\hat{y})=J(\tau_{\hat{y}}^{\lambda^{N^*}}).  
	\end{align*}
	This yields from Lemma \ref{lem:VM-V-convergence} that
	\begin{align}\label{eq:vaefcntrp0h}
		V-J(\tau^{\lambda^{N^*}}_{\hat{y}})\leq {\cal V} - \mathcal{J}(\lambda^{N^*})= {\cal V}-{\cal V}^M+{\cal V}^M-{\cal J}(\lambda^{N^*})\leq\hat{\varepsilon}_M:=\varepsilon_M+{\cal E}_M^{{N^*}}.  
	\end{align}
	Here, $\varepsilon_M=C(M^{-2}+{M}^{-1}+M^{-\frac{1}{2}})$ is provided in Lemma \ref{lem:VM-V-convergence}. This shows that $\tau^{\lambda^{N^*}}_{\hat{y}}$ is a Brownian $\hat{\varepsilon}_M$-optimal stopping time (the terminal iteration number $N^*$ can be chosen to depend on $M>1$).

	In what follows, we want to learn $\hat{y}$ via gradient ascent algorithm. To do it, we need to compute the gradient of $y\to\mathcal{Q}^{\lambda^{N^*}}(y)$, which is given in the following lemma.
	\begin{lemma}[The Gradient of $y\to\mathcal{Q}^{\lambda^{N^*}}(y)$]\label{lem:learnQy}
		Assume additionally that the terminal reward $g\in C^{1,2}([0,T]\times\R^d)$. Then, the following gradient formula holds, for any $y>0$ and $\lambda\in\U^{M,\rm m}[0,T]$,
		\begin{align}\label{eq:derivaQyy}
			\nabla_y{\cal Q}^{\lambda}(y)=\E\left[\mathbf{1}_{\{\int_0^T\lambda(s,X_s)\d s>y\}}\frac{\mu(\tau_y^{\lambda},X_{\tau_y^{\lambda}})}{\lambda(\tau_y^{\lambda},X_{\tau_y^{\lambda}})}\right],
		\end{align}
		where, the $\Fb^1$-stopping time $\tau_y^{\lambda}=\inf\{t\geq0;\int_0^t\lambda(s,X_s)\d s\geq y\}\wedge T$ for $y>0$, and the function $\mu(t,x):=f(t,x)+{\cal L}^{b,\sigma}g(t,x)$ for $(t,x)\in[0,T]\times\R^d$ with the operator ${\cal L}^{b,\sigma}\psi(t,x):=\frac{\partial}{\partial t}\psi(t,x)+b(t,x)^{\top}\nabla_x\psi(t,x)+\frac{1}{2}{\rm tr}[(\sigma\sigma^{\top})(t,x)\nabla_x^2\psi(t,x)]$ for $\psi\in C^{1,2}([0,T]\times\R^d)$.
	\end{lemma}
	
	\begin{remark}
		The assumption $g\in C^{1,2}([0,T]\times\R^d)$ can be relaxed.
		In fact, it suffices to assume that there exists a measurable function $\mu:[0,T]\times\R^d\to\R$ such that $\|\mu(\cdot,X_{\cdot})\|_{L^1}<\infty$ and the process $(\Phi_t-\int_0^t\mu(s,X_s)\d s)_{t\in[0,T]}$ is a $(\Pb^1,\Fb^1)$-martingale. Under these conditions, the gradient
		formula \eqref{eq:derivaQyy} remains valid with this choice of $\mu(t,x)$.
	\end{remark}
	
	However, the gradient formula \eqref{eq:derivaQyy} can not be directly applied since the function $\mu(t,x)=f(t,x)+{\cal L}^{b,\sigma}g(t,x)$ is unobservable ($b,\sigma$ are unknown). We therefore first seek to learn the unknown function $\mu(t,x)$ from the observable trajectory data. Recall that the reward process $\Phi=(\Phi_t)_{t\in[0,T]}$
	is observable. Thanks to the It\^o's representation \eqref{eq:dPhitmart} of the reward process in Section \ref{sec:proof_main}, the process $\Phi_t-\int_0^t\mu(s,X_s)\d s$ for $t\in[0,T]$ is an $\Fb^1$-martingale. This  observation motivates us to learn $\mu(t,x)$ by using the martingale characterization following the approach of Jia and Zhou \cite{JiaZhou22a}. More precisely, let $\mu^{\psi}(t,x)$ with $(t,x)\in[0,T]\times\R^d$ be a parameterized function with learning parameter $\psi\in\R^p$. We construct the following parameterized process given by
	\begin{align}\label{eq:Mpsi}
		{\cal M}_t^{\psi}:=\Phi_t-\int_0^t\mu^{\psi}(s,X_s)\d s,\quad \forall t\in[0,T].
	\end{align}
	By martingale orthogonality conditions, the parameterized process ${\cal M}^{\psi}=({\cal M}_t^{\psi})_{t\in[0,T]}$ is an  $\Fb^1$-martingale iff the martingale orthogonality equation $\E[\int_0^T\eta_t\d{\cal M}_t^{\psi}]=0$ holds for any $\Fb^1$-adapted continuous process $\eta=(\eta_t)_{t\in[0,T]}$ satisfying $\E[\int_0^T|\eta_t|^2\d\langle{\cal M}^{\psi}\rangle_t]<\infty$. Consequently, learning $\mu(t,x)$ can be formulated as finding an ``optimal" parameter $\psi^*\in\R^p$ such that the above martingale orthogonality conditions are satisfied. 
	
	In practice, we implement this characterization through stochastic approximation. In particular, we choose the test process in the form $\eta_s^\psi=\nabla_\psi\mu^\psi(s,X_s)$ for $s\in[0,T]$. Accordingly, for a learning rate $\zeta>0$, the parameter $\psi$ is updated according to
	\begin{align}\label{eq:SA}
		\psi \leftarrow \psi+\zeta \int_0^{T} \nabla_\psi\mu^\psi(s,X_s) \d{\cal M}_s^{\psi},
	\end{align}
	After obtaining the optimal parameter $\psi^*$, we can learn the maximizer of $y\to\mathcal{Q}^{\lambda^{N^*}}(y)$ via gradient ascent algorithm (see Lemma~\ref{lem:learnQy}). In other words, for a learning rate $\kappa>0$,
	\begin{align}\label{eq:gradient-ascent-haty}
		y\leftarrow y+\kappa \E\left[\mathbf{1}_{\{\int_0^T\lambda^{N^*}(s,X_s)\d s>y\}}\frac{\mu^{\psi^*}\left(\tau_y^{\lambda^{N^*}},X_{\tau_y^{\lambda^{N^*}}}\right)}{\lambda^{N^*}\left(\tau_y^{\lambda^{N^*}},X_{\tau_y^{\lambda^{N^*}}}\right)}\right].
	\end{align}
	
	\begin{remark}
		We can also learn the stopping gap $\Delta^{\lambda^n}(t,x)$ at the $n$-th step by using the martingale characterization approach (Jia and Zhou \cite{JiaZhou22a}). More precisely, we can implement the following steps based on the martingale characterization of $\Delta^{\lambda}(t,X_t)$. However, this approach has two limitations in the present setting. Firstly, it may require stronger regularity of the stopping gap, namely  $\Delta^{\lambda}(\cdot,\cdot)\in C^{1,2}([0,T]\times\mathbb R^d;\mathbb R)$, whereas the SMP-based approach only requires the corresponding stopping gap to be integrable. Secondly, under the current assumptions on the model coefficients, we are unable to establish a learning error estimate analogous to $S_M^{\lambda^n}(\theta)$. Consequently, the theoretical convergence guarantees and error analysis developed for our SMP-based learning approach are not available for the martingale characterization approach.
	\end{remark}
	
	We next summarize our algorithm as follows.
	
	\begin{breakablealgorithm}{SMP-Based Policy Iteration for Optimal Stopping}
		\label{alg:optimal_stopping}
		\begin{algorithmic}[1]
			\REQUIRE Initial critic parameter $\theta_0\in\R^q,\psi_0\in\R^p$, initial threshold $y_0>0$, intensity bound $M>1$, learning rate $\eta,\zeta,\kappa>0$, hyperparameters $\epsilon,\delta,\rho$, number of iterations $N_{\text{outer}},N_{\text{inner}},N_{\mu},N_{y}$.
			
			\ENSURE A near-optimal Brownian stopping time $\tau^*$.
			
			\STATE Initialize the intensity $\lambda^0(t,x)=M\boldsymbol{1}_{\{h^{\theta_0}(t,x)\geq 0\}}+M^{-1}\boldsymbol{1}_{\{h^{\theta_0}(t,x)< 0\}}$.
			\FOR{$n=0,1,\ldots,N_{\text{outer}}-1$}
			\STATE Simulate sample paths $\{X^i_{t_j}; j=1,\dots,N_T\}_{i=1}^{K}$ of the state process.
			\STATE Compute the survival weights and value of intensity control problem, for $i=1,\dots,K$,
			\begin{align*}
				L_{t_j}^{\lambda^n,i}&=\exp\left(-\sum_{\ell=1}^j\lambda^n(t_{\ell},X_{t_{\ell}}^i)(t_{\ell}-t_{\ell-1})\right),\quad j=1,\dots,N_T,\\
				\mathcal{J}^{\lambda^n}_i&=\sum_{j=1}^{N_T}L_{t_j}^{\lambda^n,i}\left(f(t_j,X_{t_j}^i)+\lambda^n(t_j,X_{t_j}^i)g(t_j,X_{t_j}^i)\right)(t_j-t_{j-1})+L_{t_{N_T}}^{\lambda^n,i}g(t_{N_T},X_{t_{N_T}}^i).
			\end{align*}
			
			\FOR{$k=1,2,\ldots,N_{\text{inner}}$}
			\STATE Compute the perturbed survival weights and the perturbed value of intensity control problem, for $l=1,\dots,q$ and $i=1,\dots, K$,
			\begin{align*}
				L_{t_j}^{\lambda^n,\theta^k,l,i}&=\exp\left(-\sum_{\ell=1}^j\left(\lambda^n(t_{\ell},X_{t_{\ell}}^i)+\epsilon\pa_{\theta_l}h^{\theta^k}(t_{\ell},X_{t_{\ell}}^i)\right)(t_{\ell}-t_{\ell-1})\right),\quad j=1,\dots,N_T.\\
				\mathcal{J}_i^{\lambda^n+\epsilon\pa_{\theta_l}h^{\theta^k}(\cdot,X_{\cdot})}&=L_{t_{N_T}}^{\lambda^n,\theta^k,l,i}g(t_{N_T},X_{t_{N_T}})+\sum_{j=1}^{N_T}L_{t_j}^{\lambda^n,\theta^k,l,i}\Big(f(t_j,X_{t_j}^i)\nonumber\\
				&\quad+(\lambda^n(t_j,X_{t_j}^i)+\epsilon \pa_{\theta_l}h^{\theta^k}(t_j,X_{t_j}^i))g(t_j,X_{t_j}^i)\Big)(t_j-t_{j-1}).
			\end{align*}
			
			\STATE Estimate the gradient of the critic loss using the SMP-based finite-difference scheme, for $l=1,\dots,q$,
			\begin{align*}
				G_l^{n,k}=\frac1K\sum_{i=1}^K\left(\sum_{j=1}^{N_T}\left(L_{t_j}^{\lambda^n,i}\pa_{\theta_l}h^{\theta^k}(t_j,X_{t_j}^i)h^{\theta^k}(t_j,X_{t_j}^i)(t_j-t_{j-1})\right)-\frac{\mathcal{J}_i^{\lambda^n+\epsilon\pa_{\theta_l}h^{\theta^k}(\cdot,X_{\cdot})}-\mathcal{J}_i^{\lambda^n}}{\epsilon}\right).
			\end{align*}
			\STATE Set the vector $G^{n,k}=(G_1^{n,k},\dots,G_q^{n,k})$ and update the critic parameter:
			\begin{align*}
				\theta^{k+1}
				=\theta^k-\eta_kG^{n,k}.   
			\end{align*}
			\ENDFOR
			\STATE The inner-loop terminates with the critic parameter $\theta^{*,n+1}$.
			\STATE Update the intensity with step $\rho_n$ as follows:
			\begin{align*}
				\lambda^{n+1}(t,x)=(1-\rho_n)\lambda^n(t,x)+\rho_n\left(M\boldsymbol{1}_{\{h^{\theta^*,n+1}(t,x)\geq 0\}}+M^{-1}\boldsymbol{1}_{\{h^{\theta^*,n+1}(t,x)<0\}}\right). 
			\end{align*}
			\ENDFOR
			\STATE The outer-loop terminates with policy $\lambda^{N^*}(t,x)$.
			\FOR{$u=0,\dots,N_{\mu}-1$}
			\STATE Simulate a sample path $\{X_{t_j}; j=1,\dots,N_T\}$ and update the parameter $\psi$ by
			\begin{align*}
				\psi^{u+1}&=\psi^u+\zeta_u\sum_{j=0}^{N_T-1}\nabla_{\psi}\mu^{\psi^u}(t_j,X_{t_j})\Big((f(t_j,X_{t_j})-\mu^{\psi^u}(t_{j},X_{t_j}))(t_{j+1}-t_j)\\
				&\quad+g(t_{j+1},X_{t_{j+1}})-g(t_j,X_{t_j})\Big).
			\end{align*}
			\ENDFOR
			\STATE The gradient descent loop terminates with parameter $\psi^*$.
			\FOR{$p=0,\ldots,N_y-1$}
			\STATE Simulate sample paths $\{X^i_{t_j}; j=1,\dots,N_T\}_{i=1}^{K}$ of the state process.
			
			\STATE Compute the stopping time $\tau_{y_p}^{\lambda^{N^*},i}$ for $i=1,\dots,K$,
			\begin{align*}
				\tau_{y_p}^{\lambda^{N^*},i}=\inf\left\{t\geq 0;~\int_0^t\lambda^{N^*}(s,X_s^i)\d s\geq y_p\right\}\wedge T.
			\end{align*}
			
			\STATE Compute the derivative $y\to\mathcal{Q}^{\lambda^{N^*}}(y)$ at $y_p$ in the way that
			\begin{align*}
				\widehat{Q}'(y_p)=\frac1K\sum_{i=1}^K\left(\boldsymbol{1}_{\left\{y_p<\int_0^T\lambda^{N^*}(s,X_s^i)\d s\right\}}\frac{\mu^{\psi^*}\left(\tau_{y_p}^{\lambda^{N^*},i},X_{\tau_{y_p}^{\lambda^{N^*},i}}^i\right)}{\lambda^{N^*}\left(\tau_{y_p}^{\lambda^{N^*},i},X_{\tau_{y_p}^{\lambda^{N^*},i}}^i\right)}\right).    
			\end{align*}
			
			\STATE Update the threshold as follows:
			\begin{align*}
				y_{p+1}=y_p+\kappa_p\widehat{Q}'(y_p).    
			\end{align*}
			\ENDFOR
			\STATE The gradient ascent loop terminates with threshold $y^*$.
			\STATE Define the stopping time:
			\begin{align*}
				\tau^*=\inf\left\{t\geq 0;~\int_0^t\lambda^{N^*}(s,X_s)\d s\geq y^*\right\}\wedge T.   
			\end{align*}
			\RETURN $\tau^*$.
		\end{algorithmic}
	\end{breakablealgorithm}

	\section{Numerical Analysis}\label{sec:numerical}
	In this section, we conduct numerical experiments to illustrate the effectiveness of Algorithm \ref{alg:optimal_stopping} for solving optimal stopping problems. 
	
	\subsection{American put option pricing}
	In this subsection, we present simulation experiments in which an RL agent learns how to exercise finite-horizon American put options with $d=r=1$. Consider a finite-horizon American put option with maturity $T>0$ and strike price $P>0$. Under the risk-neutral measure, the underlying asset price process $X=(X_t)_{t\in[0,T]}$ satisfies the following Black-Scholes model, for $t\in[0,T]$,
	\begin{align*}
		\d X_t=\gamma X_t\d t+\sigma X_t\d W_t,\quad X_0=x>0,
	\end{align*}
	where, $\gamma>0$ is the risk-free interest rate level and $\sigma>0$ is the volatility. The holder of the American put option can exercise the option at any stopping time $\tau\in{\cal T}$. The corresponding payoff function at the exercise time is given by $g(t,x)=e^{-\gamma t}(P-x)^+$. Consequently, the optimal exercise problem of this American option can be formulated as the following continuous time optimal stopping problem:
	\begin{align}\label{eq:american-put}
		V(x)=\sup_{\tau\in\mathcal T}
		\mathbb E\left[g(\tau,X_{\tau})\right]=\sup_{\tau\in\mathcal T}
		\mathbb E\left[e^{-\gamma\tau}(P-X_\tau)^+\right].
	\end{align}
	Here, the volatility $\sigma$ is unknown to the agent. To apply Lemma \ref{lem:learnQy}, since the payoff function $g(t,x)$ does not satisfy the $C^{1,2}$-regularity assumption, we propose an approximation of the original payoff by a family of smooth functions. More precisely, for any $\varepsilon>0$, let us define
	\begin{align}\label{eq:softplus}
		g_\varepsilon(t,x):=e^{-\gamma t}\varepsilon
		\log\left(1+\exp\left(\frac{P-x}{\varepsilon}\right)\right),\quad \forall x\in\R.
	\end{align}
	For every fixed $\varepsilon>0$, we have $g_\varepsilon\in C^{1,2}([0,T]\times\mathbb R)$, and it also holds that
	\begin{align}\label{eq:softplus-error}
		\|g_\varepsilon-g\|_{\infty}:=\sup_{(t,x)\in[0,T]\times\R}|g_\varepsilon(t,x)-g(t,x)|\leq \varepsilon\log 2.
	\end{align}
	This yields that the sequence of functions $g_{\varepsilon}$ converges uniformly to the original payoff $g$ on $[0,T]\times\R$ as $\varepsilon\to0$. For any fixed $\varepsilon>0$, we next consider the
	corresponding approximating optimal stopping problem given by
	\begin{align}\label{eq:american-put-epsilon}
		V_\varepsilon(x)=\sup_{\tau\in\mathcal T}
		\mathbb E\left[g_{\varepsilon}(\tau,X_{\tau})\right]=\sup_{\tau\in\mathcal T}
		\mathbb E\left[e^{-\gamma\tau}\varepsilon\log\left(
		1+\exp\left(\frac{P-X_\tau}{\varepsilon}\right)
		\right)\right].
	\end{align}
	We immediately obtain from \eqref{eq:softplus-error} that  $\|V_\varepsilon-V\|_{\infty} \leq \varepsilon\log 2$. Thus, the smooth approximation introduces a uniformly controlled error in the value function, which vanishes as $\varepsilon\to0$.
	
	We then apply Algorithm~\ref{alg:optimal_stopping} to learn the optimal Brownian stopping time. The simulation parameters are set as follows. The time horizon is $T=1.0$  with $N_T=100$ time steps and $\Delta t=T/N_T=0.01$. We set the  strike price $P=100$, the risk-free interest rate $\gamma=0.05$,  the volatility $\sigma=0.20$, the truncation parameter $M=100$, and the payoff smoothing parameter  $\varepsilon=2.0$. The initial stock price is uniformly sampled from  the interval $[70,130]$.  We parameterize the stopping-gap function by  $h^\theta(t,x)=(1-T/t)\hat{h}^\theta(t,x)$ and the auxiliary function $\mu(t,x)$ by $\mu^\psi(t,x)$, where $\hat{h}^\theta(t,x)$ and $\mu^\psi(t,x)$ are implemented by  multilayer perceptrons (MLPs)  with a single hidden layer and ReLU activations, with hidden dimensions  16 and 32, respectively. 
	
	For the SMP-based policy iteration, we use a batch size of $K=128$ trajectories. The algorithm runs for $N_{\mathrm{outer}}=10$ outer iterations, with $N_{\mathrm{inner}}=50$ inner critic updates per outer iteration. The auxiliary function $\mu^\psi$ is trained for $N_\mu=1500$ iterations and the threshold $y$ is then optimized via stochastic gradient ascent for $N_y=2000$ iterations, starting from $y_0=1.0$. For $k,n,u,p\geq 1$, the learning rates are given by
	\begin{align*}
		\eta_k &= \frac{0.003}{\sqrt{1+0.02(k-1)}}, \quad
		\rho_n = \frac{0.20}{n^{0.60}}, \quad
		\zeta_u = \frac{0.003}{u^{0.20}}, \quad
		\kappa_p =10.
	\end{align*}
	
\begin{figure}[h]
		\centering
		\includegraphics[width=0.56\textwidth]{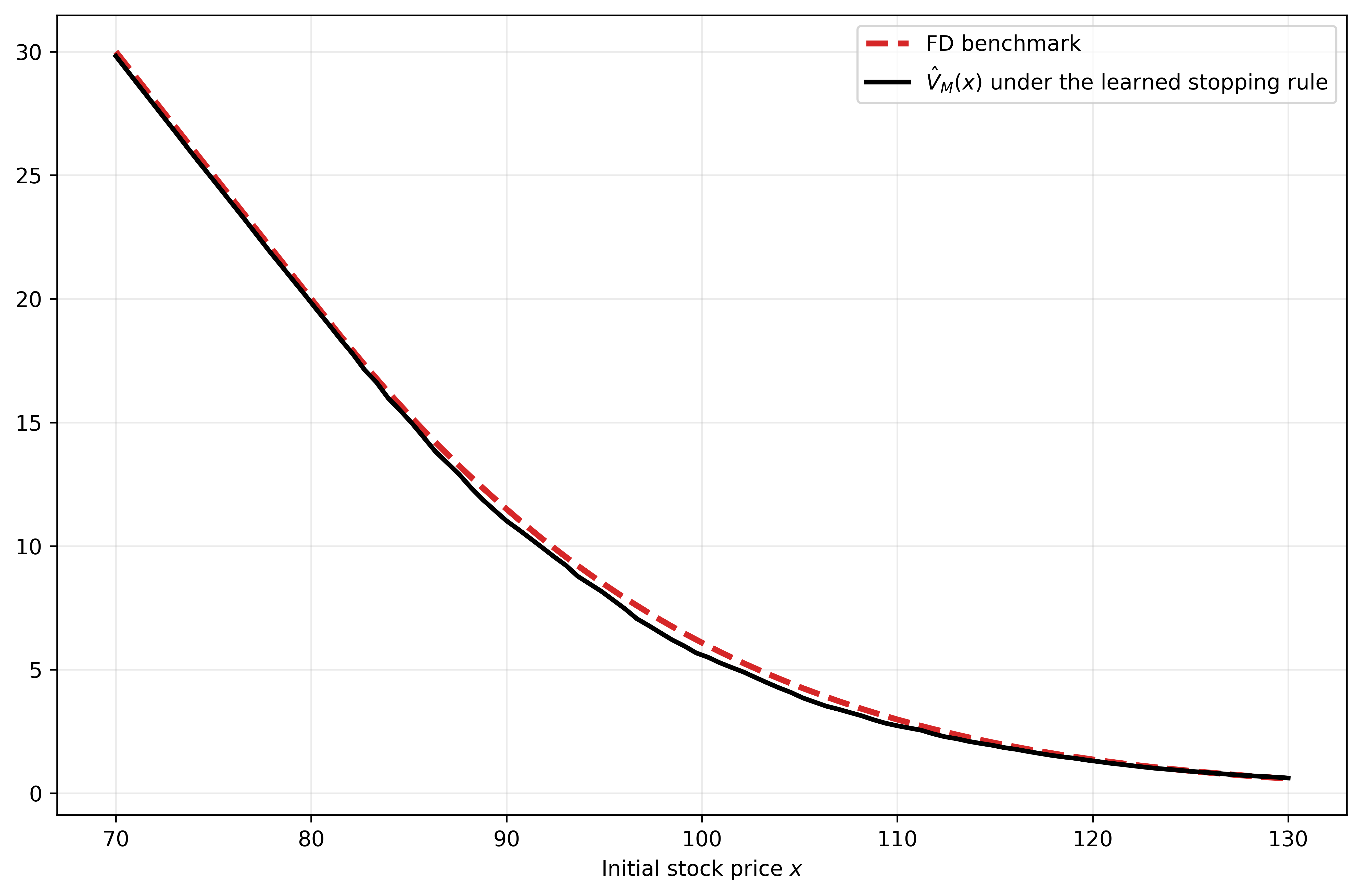}
		\caption{Value functions under the learned stopping rule vs. FD benchmark with $M=100$. }
		\label{fig:comparison-value} 
	\end{figure}

	To assess the performance of the learned stopping time in the original  optimal stopping problem \eqref{eq:american-put}, we estimate the value  function associated with the learned stopping rule using $N_{\mathrm{eval}}=1024$ independent Monte Carlo trajectories:
	\begin{align*}
		\hat{V}_M(x)=\frac{1}{N_{\mathrm{eval}}}
		\sum_{i=1}^{N_{\mathrm{eval}}}e^{-\gamma \tau_M^{*,i}}\left(P-X_{\tau_M^{*,i}}^{i}\right)^+.
	\end{align*}
	The numerical benchmark for problem \eqref{eq:american-put} is obtained using  the finite difference (FD) method. Figure~\ref{fig:comparison-value} compares 
	the estimated value function $\hat{V}_M(x)$ with the numerical benchmark value function denoted by $V_{\rm FD}(x)$. The two value functions are in close agreement, indicating that the learned stopping time achieves near-optimal performance in the original optimal stopping problem. These numerical results demonstrate the effectiveness of the proposed algorithm and provide empirical support for our  theoretical findings.
	
\begin{figure}[h]
		\centering
		\includegraphics[width=0.58\textwidth]{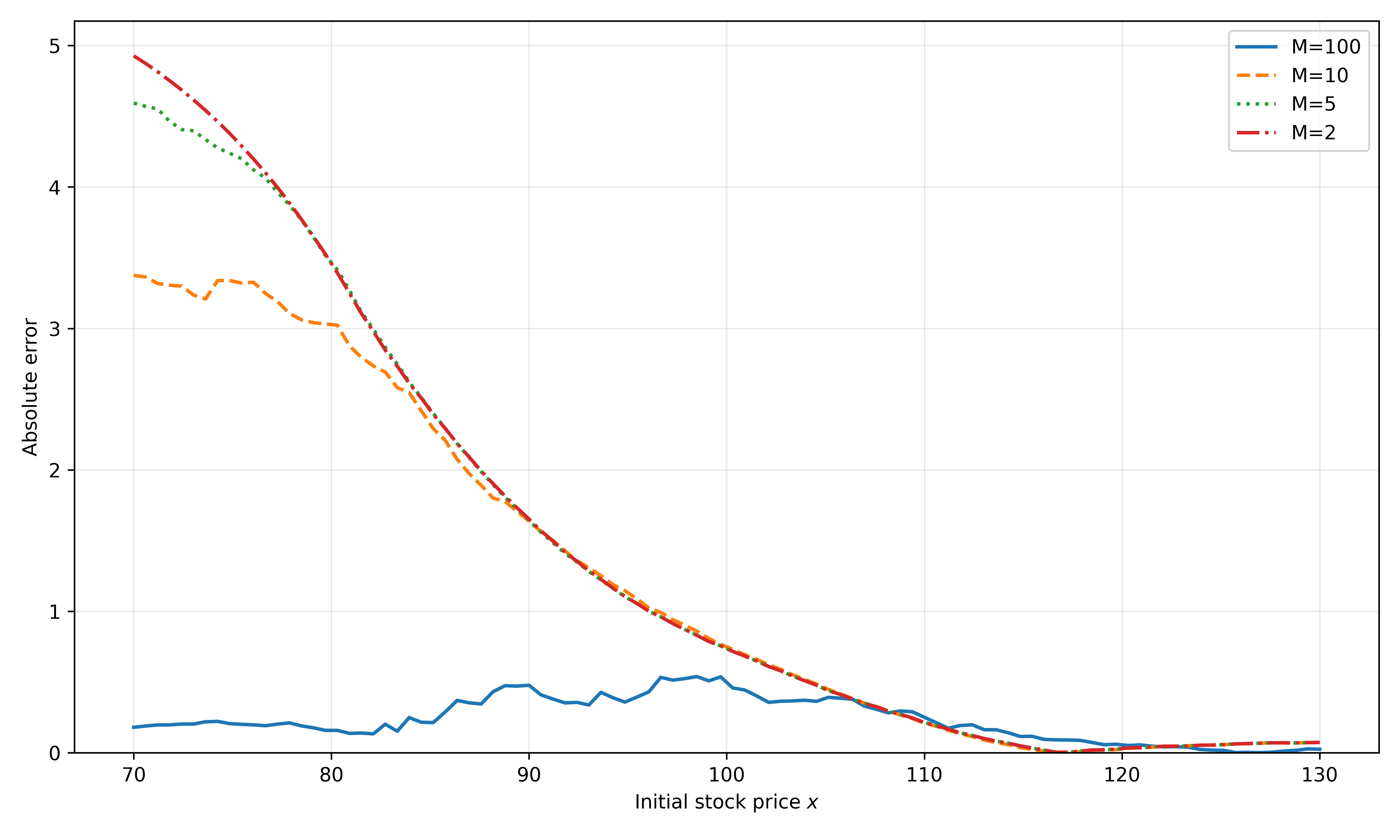}
		\caption{Comparison of absolute errors
			$|\hat{V}_M(x)-V_{\mathrm{FD}}(x)|$ for different truncation levels $M>1$. }
		\label{fig:american-put-error}
	\end{figure}

	To further investigate the effect of the truncation parameter $M>1$, we compare the value functions associated with the learned stopping  rules for different values of $M$.  Figure~\ref{fig:american-put-error} illustrates the absolute errors
	$|\hat{V}_M(x)-V_{\mathrm{FD}}(x)|$ for different truncation levels $M$. Overall, the absolute errors decrease as $M$ increases, indicating that the learned value functions progressively approach the value function of the original optimal stopping problem. This observation is consistent with the theoretical convergence result as $M\to\infty$ established in Lemma~\ref{lem:VM-V-convergence} and Theorem~\ref{thm:policy_convergence}. Although the error bound in Theorem~\ref{prop:relaxed-policy-convergence} depends on $M$, our numerical results suggest that the value functions associated with the learned stopping times approximate the FD benchmark well, with no apparent deterioration in accuracy as $M$ increases.
	
	\subsection{First change-point detection}
	In this subsection, we consider the five-dimensional first change-point
	detection problem, adapted from Ekstr{\"o}m and Wang \cite{Ekstrom22}, to further investigate the effectiveness of our algorithm in solving high-dimensional optimal stopping problems with $d=r=5$.
	
	For $i=1,\ldots,d$, let $X_t^i$ denote the posterior probability that the $i$-th
	change-point has occurred by time $t\in[0,T]$. The posterior probability process $X^i=(X_t^i)_{t\in[0,T]}$ satisfies that
	\begin{align*}
		\d X_t^i=b_i(1-X_t^i)\d t+\sigma_iX_t^i(1-X_t^i)\d W_t^i,
		\quad X_0^i=x_i\in[0,1],
	\end{align*}
	where $b_i>0$ is the change-point arrival rate, $\sigma_i>0$ is the signal strength and $W^i=(W_t^{i})_{t\in[0,T]}$ for $i=1,\ldots,d$ are independent scalar Brownian motions. Denote by $X_t=(X_t^1,\ldots,X_t^d)^{\top}$ for $t\in[0,T]$.
	
	Let us introduce that $q(x):=1-\prod_{i=1}^d(1-x_i)$ for $x\in[0,1]^d$. Then, $q(X_t)$ represents the conditional probability that at least one change-point has occurred by time $t\in[0,T]$. The objective is to detect the first change-point as quickly as possible while minimizing the probability of a false alarm. Accordingly, the first change-point detection problem can be formulated as the following optimal stopping problem, for $x\in[0,1]^d$,
	\begin{align*}
		C(x):=\inf_{\tau\in\mathcal T}\mathbb E\left[1-q(X_\tau)+c\int_0^\tau q(X_s)\d s\right].
	\end{align*}
	Here, $c>0$ denotes the penalty per unit of detection delay. To conform with the maximization framework adopted in this paper,
	define $V(x):=1-C(x)$, and consider the following equivalent optimal stopping problem, for $x\in[0,1]^d$,
	\begin{align}\label{eq:change-point}
		V(x)=\sup_{\tau\in\mathcal T}\mathbb E\left[q(X_\tau)-c\int_0^\tau q(X_s)\d s\right].
	\end{align}
	The corresponding running and terminal reward functions are
	$f(t,x)=-cq(x)$ and $g(t,x)=q(x)$, respectively. Here, both $b=(b_1,\ldots,b_d)^{\top}$ and $\sigma=(\sigma_1,\ldots,\sigma_d)^{\top}$ are unknown.
	
	We then apply Algorithm~\ref{alg:optimal_stopping} to learn the
	optimal Brownian stopping time. The simulation parameters are
	set as follows. The time horizon is $T=2.0$, the detection delay
	penalty is $c=1.0$, and we set $b_i=0.2$ and $\sigma_i=1.0$
	for all $i=1,\ldots,d$. We consider different truncation levels
	$M\in\{5,10,20,50,100\}$. To control the time discretization error
	associated with the intensity process, we set the number of time
	steps to $N_T=\max\{400,\lceil40M\rceil\}$.
	
	We parameterize the stopping-gap function by
	$h^\theta(t,x)=(1-t/T)\hat h^\theta(t,x)$ and the auxiliary
	function $\mu(t,x)$ by $\mu^\psi(t,x)$, where $\hat h^\theta$ and $\mu^\psi$ are implemented by MLPs with a single hidden layer and $\tanh$ activations, with hidden dimensions 6 and 24, respectively. For the SMP-based policy iteration, we use a batch size of $K=32$ trajectories. The algorithm runs for $N_{\mathrm{outer}}=16$ outer iterations, with $N_{\mathrm{inner}}=12$ inner critic updates per outer iteration. The auxiliary function $\mu^\psi$ is trained for $N_\mu=3200$ iterations, and the threshold $y$ is subsequently optimized via stochastic gradient ascent for $N_y=180$ iterations, starting from $y_0=1.0$. For $k,n,u,p\geq1$, the learning rates are given by
	\begin{align*}
		\eta_k&=\frac{0.08}{\sqrt{1+0.02(k-1)}},\quad
		\rho_n=\frac{0.65}{n^{0.30}},\quad
		\zeta_u=\frac{0.02}{u^{0.20}},\quad
		\kappa_p=\frac{0.05}{p^{0.60}}.
	\end{align*}
	
\begin{figure}[h]
		\centering
		\includegraphics[width=0.58\textwidth]{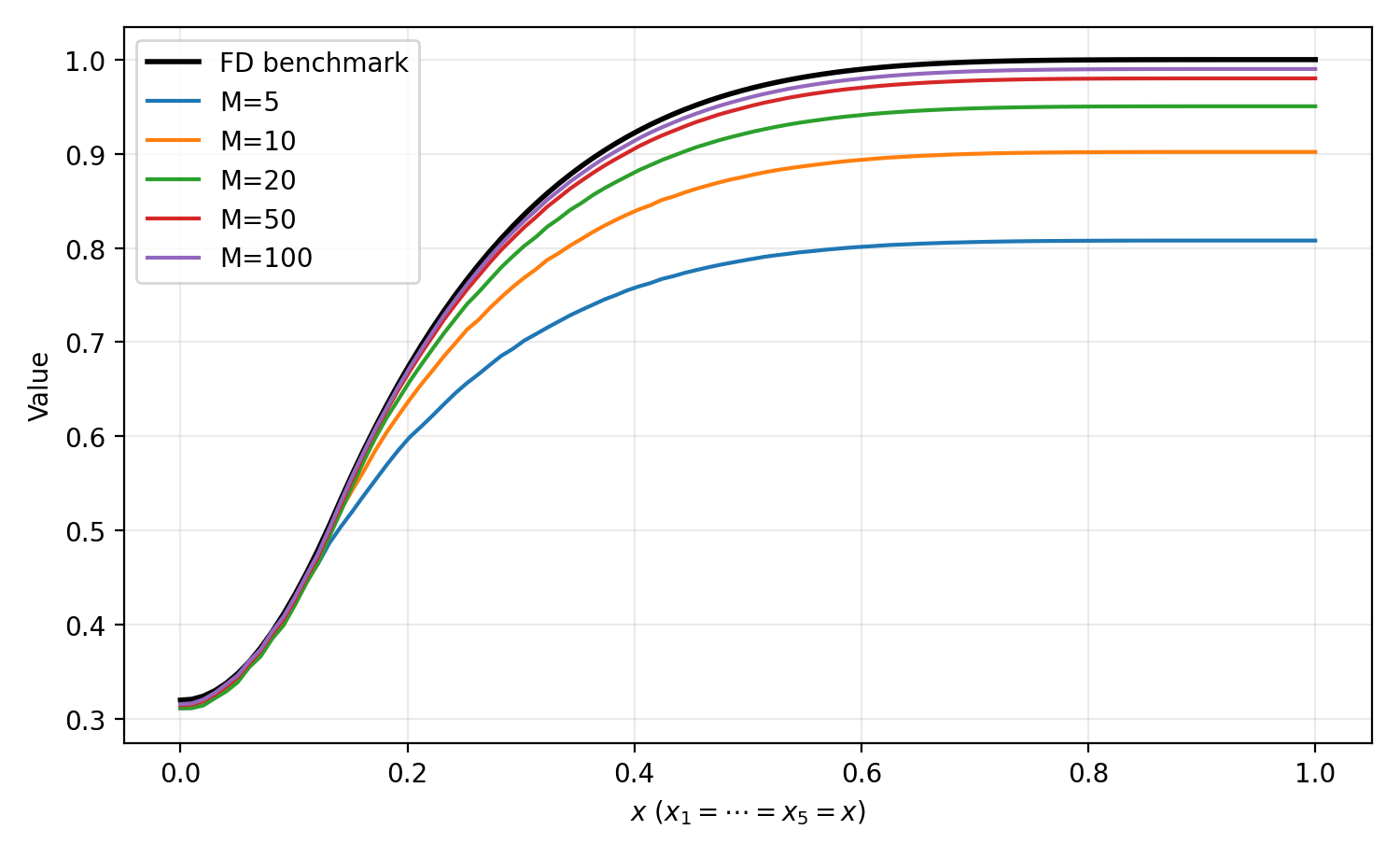}
		\caption{Comparison of value functions under learned stopping rules
			for different truncation levels $M$ with the FD benchmark. }
		\label{fig:comparison}
	\end{figure}

	To assess the performance of the learned stopping time in the
	original optimal stopping problem \eqref{eq:change-point}, we
	estimate its associated value function using
	$N_{\mathrm{eval}}=4096$ independent Monte Carlo trajectories.
	The numerical benchmark is obtained by solving the HJB
	variational inequality associated with problem
	\eqref{eq:change-point} using the FD method. Figure~\ref{fig:comparison} compares the value functions associated with the learned stopping rules for different truncation levels \(M\) with the FD benchmark. Overall, as \(M\) increases, the learned value functions progressively approach the numerical benchmark. This observation is consistent with our theoretical convergence results and demonstrates the effectiveness of the proposed algorithm in solving high-dimensional optimal stopping problems.
	
	\section{Proofs of Results}\label{sec:proof_main}
	
	In this section, we present the proofs of lemmas, propositions and theorems appearing in the previous sections. Throughout the section, Assumption \ref{ass1} is assumed to hold.
	
	\begin{proof}[Proof of Lemma \ref{lem:valuefcnrep}]
		On one hand, we have from the tower property and the fact $\xi\sim \mathrm{Exp}(1)$ that
		\begin{align*}    
			\mathcal{J}(\lambda)=\E\left[\Phi_{\tau_{\xi}^{\lambda}}\right]=\E\left[\E\left[\Phi_{\tau_{\xi}^{\lambda}}\middle|\xi\right]\right]=\int_0^{\infty}\mathcal{Q}^{\lambda}(y)e^{-y}\d y.
		\end{align*}
		This verifies the validity of the first equality in \eqref{eq:vaefcntrp0}. On the other hand, note that $\{\tau_{\xi}^{\lambda}>s\}=\{T>s\}\cap\{\xi>\int_0^s\lambda_r\d r\}$ for $s\geq0$. It follows from the fact $\xi\perp\F_T^1$ that
		\begin{align}\label{eq:condexp1}
			{\cal J}(\lambda)&=\E\left[\int_0^{\tau_{\xi}^{\lambda}}f(s,X_s)\d s+g(\tau_{\xi}^{\lambda},X_{\tau_{\xi}^{\lambda}})\right]\nonumber\\
			&=\E\left\{\E\left[\int_0^{\infty}{\bf1}_{\{\tau_{\xi}^{\lambda}>s\}}f(s,X_s)\d s+g(\tau_{\xi}^{\lambda},X_{\tau_{\xi}^{\lambda}})\bigg|\F_T^1\right]\right\}\nonumber\\
			&=\E\left\{\E\left[\int_0^{T}{\bf1}_{\{\xi>\int_0^s\lambda_r\d r\}}f(s,X_s)\d s\bigg|\F_T^1\right]\right\}+\E\left\{\E\left[g(\tau_{\xi}^{\lambda},X_{\tau_{\xi}^{\lambda}})\bigg|\F_T^1\right]\right\}.
		\end{align}
		The first term of RHS of the equality \eqref{eq:condexp1} equals to $\E[\int_0^{T}L_s^{\lambda}f(s,X_s)\d s]$ since $\xi\sim{\rm Exp}(1)$. For the second term of RHS of the equality \eqref{eq:condexp1}, let $\hat{\tau}_{\xi}^{\lambda}:=\inf\{t\geq 0;~\int_0^t\lambda_r\d r\geq\xi\}$, we have 
		\begin{align*}
			\E\left[g(\tau_{\xi}^{\lambda},X_{\tau_{\xi}^{\lambda}})\bigg|\F_T^1\right]&=\E\left[g(\hat{\tau}_{\xi}^{\lambda},X_{\hat{\tau}_{\xi}^{\lambda}}){\bf 1}_{\{\hat{\tau}_{\xi}^{\lambda}\leq T\}}\bigg|\F_T^1\right]+\E\left[g(T,X_T){\bf 1}_{\{\hat{\tau}_{\xi}^{\lambda}> T\}}\bigg|\F_T^1\right]\nonumber\\
			&=-\E\left[\int_0^Tg(s,X_s)\d L_s^{\lambda}\bigg|\F_T^1\right]+\E\left[L_T^{\lambda}g(T,X_T)\bigg|\F_T^1\right]\nonumber\\
			&=\E\left[\int_0^TL_s^{\lambda}\lambda_sg(s,X_s)\d s+L_T^{\lambda}g(T,X_T)\bigg|\F_T^1\right].
		\end{align*}
		This finishes the verification of the representation \eqref{eq:vaefcntrp0}. 
	\end{proof}
	
	\begin{proof}[Proof of Lemma \ref{lem:BSDElambdadecoupling}]
		Consider the adjoint equation of the survival component in \eqref{eq:BSDE}, which is given by, for $t\in[0,T]$, 
		\begin{align}\label{eq:BSDE_p_int}
			\d P_t^{\lambda}=-(-\lambda_tP_t^{\lambda}+C_t)\d t+Q_t^{\lambda}\d W_t,\quad  P_T^{\lambda}=g(T,X_T),  
		\end{align}
		where, the process $C_t:=f(t,X_t)+\lambda(t,X_t)g(t,X_t)$ for $t\in[0,T]$. Thus, BSDE \eqref{eq:BSDE_p_int} is a linear BSDE. By Assumption \ref{ass1}, we have $\E[\int_0^T|C_t|^2\d t]<\infty$ and $\E[|g(T,X_T)|^2]<\infty$. This yields that this BSDE has a unique solution $(P^{\lambda},Q^{\lambda})=(P_t^{\lambda},Q_t^{\lambda})_{t\in[0,T]}$ satisfying $\E[\sup_{t\in[0,T]}|P_t^{\lambda}|^2]<\infty$ and $\E[\int_0^T|Q_t^{\lambda}|^2\d t]<\infty$. Moreover, in light of Proposition 6.2.1 of Pham \cite{Pham} and Markov property, this solution has the following closed-form given by, for $t\in[0,T]$,
		\begin{align*}
			P_t^{\lambda}&=\E\left[\exp\left(-\int_t^T\lambda_r\d r\right)g(T,X_T)+\int_t^T\exp\left(-\int_t^r\lambda_u\d u\right)C_r\d r\bigg|\F_t^1\right]=v^{\lambda}(t,X_t).
		\end{align*}  
		Thus, we complete the proof of the lemma.
	\end{proof}
	
	\begin{proof}[Proof of Lemma \ref{lem:VM-V-convergence}]
		For any $\varepsilon>0$, let $\tau_{\varepsilon}\in{\cal T}$ be an $\varepsilon$-optimal stopping time of problem \eqref{eq:reward_func}, i.e., 
		\begin{align*}
			J(\tau_{\varepsilon})=\E[\Phi_{\tau_{\varepsilon}}]=\E\left[\int_0^{\tau_{\varepsilon}} f(s,X_s)\d s+g\left(\tau_{\varepsilon},X_{\tau_{\varepsilon}}\right)\right]\geq {\cal V}-\varepsilon.    
		\end{align*}
		Introduce an admissible stopping intensity $\lambda_t^{\varepsilon}:=M\mathbf{1}_{\{t>\tau_\varepsilon\}}+M^{-1}\mathbf{1}_{\{t\leq\tau_\varepsilon\}}$ for $t\in[0,T]$. Then, $\lambda^{\varepsilon}=(\lambda_t^{\varepsilon})_{t\in[0,T]}\in\mathscr{U}^M[0,T]$. Hence, in terms of \eqref{eq:taulambda}, we have
		
			\begin{align}\label{eq:taulvarepsilon}
				\tau^{\lambda^{\varepsilon}}&=\inf\left\{t\geq 0;~\int_0^t\lambda_r^{\varepsilon}\d r\geq\xi\right\}\wedge T=\inf\left\{t\geq 0;~t\geq\left(1-\frac{1}{M^2}\right)\tau_{\varepsilon}+\frac{\xi}{M}\right\}\wedge T\nonumber\\
				&=\hat{\tau}_{\varepsilon}\wedge T:=\left[\left(1-\frac{1}{M^2}\right)\tau_{\varepsilon}+\frac{\xi}{M}\right]\wedge T.   
		\end{align}
		This yields that
\begin{align*}
	{\cal V}^M \geq{\cal J}(\lambda^{\varepsilon})
	=\E\left[\int_0^{\tau^{\lambda^{\varepsilon}}} f(s,X_s)\d s
	+g\left(\tau^{\lambda^{\varepsilon}},X_{\tau^{\lambda^{\varepsilon}}}\right)\right].
\end{align*}
In the following, let $C>0$ be a generic constant independent of $M$, which may be different from line to line. As a consequence,
\begin{align}\label{eq:estimate01}
	0&\leq {\cal V}-{\cal V}^{M}\nonumber\\
	&\leq \E\left[\int_0^{\tau_{\varepsilon}} f(s,X_s)\d s+g(\tau_{\varepsilon},X_{\tau_{\varepsilon}})\right]
	-\E\left[\int_0^{\tau^{\lambda^{\varepsilon}}} f(s,X_s)\d s
	+g\left(\tau^{\lambda^{\varepsilon}},X_{\tau^{\lambda^{\varepsilon}}}\right)\right]+\varepsilon\nonumber\\
	&=\varepsilon+\E\left[\int_{\tau^{\lambda^{\varepsilon}}}^{\tau_{\varepsilon}} f(s,X_s)\d s\right]
	+\E\left[g\left(\tau_{\varepsilon},X_{\tau_{\varepsilon}}\right)
	-g\left(\tau^{\lambda^{\varepsilon}},X_{\tau^{\lambda^{\varepsilon}}}\right)\right]\nonumber\\
	&\leq\varepsilon+\E\left[\sup_{s\in[0,T]}\left|f(s,X_s)\right|\left|\tau^{\lambda^{\varepsilon}}-\tau_{\varepsilon}\right|\right]
	+C\,\E\left[\left|\tau^{\lambda^{\varepsilon}}-\tau_{\varepsilon}\right|\right]
	+C\,\E\left[\left|X_{\tau_{\varepsilon}}-X_{\tau^{\lambda^{\varepsilon}}}\right|\right].
\end{align}
Here the third line uses the convention for oriented integrals, i.e.,
$\int_{\tau^{\lambda^{\varepsilon}}}^{\tau_{\varepsilon}} = -\int_{\tau_{\varepsilon}}^{\tau^{\lambda^{\varepsilon}}}$,
so no ordering of $\tau^{\lambda^{\varepsilon}}$ and $\tau_{\varepsilon}$ is needed; the fourth line follows from
$\E[\int_a^b f]\le \E|\int_a^b f|\le \E[\sup_{s\in[0,T]}|f(s,X_s)|\,|b-a|]$
and from the Lipschitz continuity of $g$ in $(t,x)$.
		
		Using the squared growth of the running reward function $x\to f(t,x)$, one has from Assumption \ref{ass1} that (c.f. Karatzas and Shreve \cite{KaratzasShreve88})
		\begin{align*}
			\E\left[\sup_{s\in[0,T]}\left|f(s,X_s)\right|\right]\leq C\left\{1+\E\left[\sup_{s\in[0,T]}\left|X_s\right|^2\right]\right\}\leq C\left\{1+\E\left[\left|\eta\right|^2\right]\right\}.    
		\end{align*}
		On the other hand, we have from \eqref{eq:taulvarepsilon} that
		\begin{align}\label{eq:lengthtime}
				\left|\tau^{\lambda^{\varepsilon}}-\tau_{\varepsilon}\right|&=\left|\hat{\tau}_{\varepsilon}\wedge T-\tau_{\varepsilon}\wedge T\right|\leq \left|\hat{\tau}_{\varepsilon}-\tau_{\varepsilon}\right|\wedge T=\left|\frac{\xi}{M}-\frac{\tau_{\varepsilon}}{M^2}\right|\wedge T\leq\left|\frac{\xi}{M}-\frac{\tau_{\varepsilon}}{M^2}\right|\nonumber\\
				&\leq \frac{\xi}{M}+\frac{T}{M^2}.
			\end{align}
			This yields that $\E[|\tau^{\lambda^{\varepsilon}}-\tau_{\varepsilon}|]\leq\frac{1}{M}+\frac{T}{M^2}$. It follows from the fact $\xi\perp\F_T^1$ and \eqref{eq:lengthtime} that 
		\begin{align*}
			&\E\left[\sup_{s\in[0,T]}\left|f(s,X_s)\right|\left|\tau^{\lambda^{\varepsilon}}-\tau_{\varepsilon}\right|\right]=\E\left[\sup_{s\in[0,T]}\left|f(s,X_s)\right|\left(\frac{\xi}{M}+\frac{T}{M^2}\right)\right]\nonumber\\
			&\qquad\leq\left(\frac{1}{M}+\frac{T}{M^2}\right)\E\left[\sup_{s\in[0,T]}\left|f(s,X_s)\right|\right]\leq C\left\{1+\E\left[\left|\eta\right|^2\right]\right\}\left(\frac{1}{M}+ \frac{1}{M^2}\right).    
		\end{align*}
		Lastly, we estimate $\E[|X_{\tau_{\varepsilon}}-X_{\tau^{\lambda^{\varepsilon}}}|]$. In fact, we have from \eqref{eq:SDE_dynamic} that
		\begin{align*}
			\E\left[\left|X_{\tau_{\varepsilon}}-X_{\tau^{\lambda^{\varepsilon}}}\right|\right]\leq \E\left[\int_{\tau_{\varepsilon}}^{\tau^{\lambda^{\varepsilon}}}\left|b(s,X_s)\right|\d s\right]+\E\left[\left|\int_{\tau_{\varepsilon}}^{\tau^{\lambda^{\varepsilon}}}\sigma(s,X_s)\d W_s\right|\right].
		\end{align*}
		For the first term of RHS of the above display, we have from \eqref{eq:lengthtime} and the linear growth of $x\to b(t,x)$ that
		\begin{align*}
			&\E\left[\int_{\tau_{\varepsilon}}^{\tau^{\lambda^{\varepsilon}}}\left|b(s,X_s)\right|\d s\right]=\E\left[\int_{0}^T \left|b(s,X_s)\right|\mathbf{1}_{\{\tau_{\varepsilon}<s\leq\tau^{\lambda^{\varepsilon}}\}}\d s\right]\leq\E\left[\sup_{s\in[0,T]}\left|b(s,X_s)\right|\left|\tau^{\lambda^{\varepsilon}}-\tau_{\varepsilon}\right|\right]\nonumber\\
			&\quad \leq\E\left[\sup_{s\in[0,T]}\left|b(s,X_s)\right|\left(\frac{\xi}{M}+\frac{T}{M^2}\right)\right]\leq C\left\{1+\E\left[\left|\eta\right|\right]\right\}\left(\frac{1}{M}+\frac{T}{M^2}\right).
		\end{align*}
		Similarly, we have from the It\^o's isometry and \eqref{eq:lengthtime} that
		\begin{align*}
			&\E\left[\left|\int_{\tau_{\varepsilon}}^{\tau^{\lambda^{\varepsilon}}}\sigma(s,X_s)\d W_s\right|^2\right] =\E\left[\int_0^T\left|\sigma(s,X_s)\right|^2\mathbf{1}_{\{\tau_{\varepsilon}<s\leq\tau^{\lambda^{\varepsilon}}\}}\d s\right]\leq\E\left[\sup_{s\in[0,T]}\left|\sigma(s,X_s)\right|^2\left|\tau^{\lambda^{\varepsilon}}-\tau_{\varepsilon}\right|\right]\nonumber\\
			&\quad\leq\E\left[\sup_{s\in[0,T]}\left|\sigma(s,X_s)\right|^2\left|\tau^{\lambda^{\varepsilon}}-\tau_{\varepsilon}\right|\right]\leq C\left\{1+\E\left[\left|\eta\right|^2\right]\right\}\left(\frac{1}{M}+\frac{T}{M^2}\right)\leq \frac{C}{M}\left\{1+\E\left[\left|\eta\right|^2\right]\right\},
		\end{align*}
		where, we used the inequality $M^{-2}T\leq M^{-1}T$ for $M>1$. This concludes that
		\begin{align}\label{eq:differencestopX}
			\E\left[\left|X_{\tau_{\varepsilon}}-X_{\tau^{\lambda^{\varepsilon}}}\right|\right]\leq C\left\{\frac{1+\E\left[\left|\eta\right|\right]}{M^2}+\frac{1+\E\left[\left|\eta\right|\right]}{M}+\frac{\sqrt{1+\E\left[|\eta|^2\right]}}{\sqrt{M}}\right\}.
		\end{align}
		Combining \eqref{eq:estimate01} and \eqref{eq:differencestopX}, we arrive at, for any $\varepsilon>0$,
		\begin{align}\label{eq:estimate1}
			&0\leq {\cal V}-{\cal V}^{M}\leq C\left\{\frac{1+\E\left[|\eta|\right]+\E\left[|\eta|^2\right]}{M^2}+\frac{1+\E\left[|\eta|\right]+\E\left[|\eta|^2\right]}{M} + \frac{\sqrt{1+\E\left[|\eta|^2\right]}}{\sqrt{M}}\right\}+\varepsilon.
		\end{align}
		Then, the desired estimate \eqref{eq:estimate2} follows from the arbitrariness of $\varepsilon>0$ in \eqref{eq:estimate1}.
	\end{proof}
	
	\begin{proof}[Proof of Lemma \ref{lem:cappedintensity}]
		Consider the adjoint equation of the survival component in \eqref{eq:BSDE}, which is given by the following linear scalar BSDE, for any $t\in[0,T]$ and $\lambda\in\mathscr{U}^M[0,T]$,
		\begin{align}\label{eq:BSDE_p}
			P_t^{\lambda}=g(T,X_T)+\int_t^T\{f(s,X_s)+\lambda_s(g(s,X_s)-P_s^{\lambda})\}\d s-\int_t^T Q_s^{\lambda}\d W_s.
		\end{align} 
		Let $(t,x)\in[0,T]\times\R^d$. Define $a^*=M\mathbf{1}_{\{g(t,x)\geq p\}}+M^{-1}\mathbf{1}_{\{g(t,x)<p\}}$ for $p\in\R$. Then, for any $a\in[M^{-1},M]$ and $p\in\R$, 
		\begin{align}\label{eq:fgineqmax}
				f(t,x)+ a(g(t,x)-p)&\leq f(t,x)+ a^*(g(t,x)-p)\nonumber\\
				&= f(t,x)+M(g(t,x)-p)^+ - \frac{1}{M}(g(t,x)-p)^{-}.
		\end{align}
		Note that the generator of BSDE \eqref{eq:BSDEPmM} is Lipschitz w.r.t. solution variable $p\in\R$. Then, the comparison theorem of the scalar BSDE implies $P_t^{\lambda}\leq P_t^{M}$ for all $t\in[0,T]$, a.s.. Using the equality relation in \eqref{eq:fgineqmax}, we have from the uniqueness of linear BSDE that $\lambda^{M,*}=(\lambda_t^{M,*})_{t\in[0,T]}\in\mathscr{U}^M[0,T]$ is an optimal stopping intensity  of the truncated problem the truncated problem $\sup_{\lambda\in\mathscr{U}^M[0,T]}{\cal J}(\lambda)$.
	\end{proof}
	
	\begin{proof}[Proof of Proposition \ref{prop:performancediff}]
		By applying It\^o's formula to $L_t^{\phi}P_t^{\lambda}$ for $t\in[0,T]$, we have, for any $t\in[0,T]$,
		\begin{align*}
			L_T^{\phi}P_T^{\lambda}&=L_t^{\phi}P_t^{\lambda} - \int_t^T L_r^{\phi}\{f(r,X_r)+\lambda_rg(r,X_r)+P_r^{\lambda}(\phi_r-\lambda_r)\}\d r+\int_t^T L_r^{\phi}Q_r^{\lambda}\d W_r.
		\end{align*}
		Taking conditional expectation $\E[\cdot|\F_t^1]$ on both sides of the above display, we have from \eqref{eq:vlambdaJ} that, for any $t\in[0,T]$,
		\begin{align*}
			L_t^{\phi}P_t^{\lambda}&=
			\E\left[\int_t^T L_s^{\phi}\left(f(s,X_s)+\lambda_sg(s,X_s)+P_s^{\lambda}(\phi_s-\lambda_s)\right)\d s+L_T^{\phi}g(T,X_T)\Big|\F_t^1\right].
		\end{align*}
		It follows from \eqref{eq:BSDE} that, for any $t\in[0,T]$, 
		\begin{align*}
			L_t^{\phi}P_t^{\phi}=\E\left[\int_t^TL_s^{\phi}(f(s,X_s)+\phi_sg(s,X_s))\d s+L_T^{\phi}g(T,X_T)\Big|\F_t^1\right].    
		\end{align*}
		By making the difference of the above two displays, we arrive at the equality given by, for any $t\in[0,T]$,
		\begin{align}\label{eq:differenceJ03}
			P_t^{\phi}-P_t^{\lambda}&= \E\left[\int_t^T(L_t^{\phi})^{-1}L_s^{\phi}(\phi_s-\lambda_s)(g(s,X_s)-P_s^{\lambda})\d s\Big|\F_t^1\right].
		\end{align}
		Furthermore, it follows from \eqref{eq:differenceJ03} and \eqref{eq:decouplfields} that $P_t^{\lambda}=v^{\lambda}(t,X_t)$, $\d t\times\Pb$-a.s.. Plugging this into \eqref{eq:differenceJ03} to conclude the validity of \eqref{eq:differenceJ02}.
	\end{proof}
	
	\begin{proof}[Proof of Lemma \ref{lem:stoptimedelta}]
		We rewrite the linear BSDE \eqref{eq:BSDE_p} by using the stopping gap $\Delta^{\lambda}(t,x)$ as follows, for any $t\in[0,T]$ and $\lambda\in\mathscr{U}^{{\rm m}}[0,T]$,
		\begin{align*}
			P_t^{\lambda}=g(T,X_T)+\int_t^T\{f(s,X_s)+\lambda_s\Delta^{\lambda}(s,X_s)\}\d s-\int_t^T Q_s^{\lambda}\d W_s.
		\end{align*} 
		Define $\tilde{P}_t:=P_t^{\lambda}+\int_0^tf(s,X_s)\d s$ for $t\in[0,T]$. Then, it holds that
		\begin{align*}
			\d\tilde{P}_t^{\lambda}=-\lambda_t\Delta^{\lambda}(t,X_t)\d t+Q_t^{\lambda}\d W_t.
		\end{align*} 
		Note that $P_t^{\lambda}=v^{\lambda}(t,X_t)$ by using \eqref{eq:decouplfields}. We have from the optional sampling theorem and \eqref{eq:vlambdaJ} that, for any $\tau\in{\cal T}$,
		\begin{align}\label{eq:EPhitau}
			J(\tau)&=\E\left[g(\tau,X_{\tau}) + \int_0^{\tau}f(s,X_s)\d s\right] = \E\left[\tilde{P}_{\tau}^{\lambda}\right] + \E\left[\Delta^{\lambda}(\tau,X_{\tau})\right]\nonumber\\
			&=\E\left[v^{\lambda}(0,X_0)\right]+\E\left[\Delta^{\lambda}(\tau,X_{\tau})\right]-\E\left[\int_0^{\tau}\lambda_s\Delta^{\lambda}(s,X_s)\d s\right]\nonumber\\
			&={\cal J}(\lambda)+\E\left[\Delta^{\lambda}(\tau,X_{\tau})\right]-\E\left[\int_0^{\tau}\lambda_s\Delta^{\lambda}(s,X_s)\d s\right].
		\end{align}
		Recall $\tau_{\lambda}:= \inf\{t\geq0;~\Delta^{\lambda}(t,X_t)\geq0\}\wedge T$ given in \eqref{eq:stoptimetaulambda}. Taking $\tau=\tau_{\lambda}$ in \eqref{eq:EPhitau}, we have
		\begin{align*}
			V=\sup_{\tau\in{\cal T}}J(\tau)\geq J(\tau_{\lambda})={\cal J}(\lambda)+\underbrace{\E\left[\Delta^{\lambda}(\tau_{\lambda},X_{\tau_{\lambda}})\right]-\E\left[\int_0^{\tau_{\lambda}}\lambda_s\Delta^{\lambda}(s,X_s)\d s\right]}_{\geq0}\geq{\cal J}(\lambda).
		\end{align*}
		This yields that $0\leq V-J(\tau_{\lambda})\leq\delta_{\lambda}$ since ${\cal V}=V$. Thus, we complete the proof of the lemma.
	\end{proof}
	
	\begin{proof}[Proof of Corollary \ref{coro:approximatestoptime}]
		Recall $\delta_{\lambda}={\cal V}-{\cal J}(\lambda)$ in Lemma \ref{lem:stoptimedelta} for $\lambda\in\mathscr{U}^{{\rm m}}[0,T]$. Then, it follows from Lemma \ref{lem:cappedintensity} and Lemma \ref{lem:VM-V-convergence} that  
		\begin{align}\label{eq:approintensityM}
			\delta_{\lambda^{M,*}}={\cal V}-{\cal J}(\lambda^{M,*})={\cal V}-{\cal V}^M \leq \varepsilon_M.   
		\end{align}
		Thus, Lemma \ref{lem:stoptimedelta} implies that $\tau_{\lambda^{M,*}}$ is a $\varepsilon_M$-optimal stopping time of the prime problem \eqref{eq:reward_func}.
	\end{proof}
	
	\begin{proof}[Proof of Lemma \ref{lem:policy_improvement}]
		In light of Proposition~\ref{prop:performancediff}, we have, for any $\lambda\in\mathscr{U}^{M,{\rm m}}[0,T]$, 
		\begin{align}\label{eq:value_dif}
			\mathcal{J}(\hat{\lambda}^M)-\mathcal{J}(\lambda)=\E\left[\int_0^TL_s^{\hat{\lambda}^M}(\hat{\lambda}^M(s,X_s)-\lambda(s,X_s))\Delta^{\lambda}(s,X_s)\d s \right].
		\end{align}
		By construction of $\hat{\lambda}^M$, for any $s\in[0,T]$, we can deduce that,  on $\{\Delta^{\lambda}(s,X_s)\geq0\}$, it holds that $\hat{\lambda}^M_s-\lambda_s= M - \lambda_s\geq0$; while on $\{\Delta^{\lambda}(s,X_s)<0\}$, one has $\hat{\lambda}^M_s-\lambda_s= M^{-1} - \lambda_s\leq0$.
		Note that $L_s^{\hat{\lambda}^M}>0$ for all $s\in[0,T]$. Then, it follows from \eqref{eq:value_dif} that $\mathcal{J}(\hat{\lambda}^M)\geq\mathcal{J}(\lambda)$. This proves the policy improvement.
	\end{proof}
	
	\begin{proof}[Proof of Theorem \ref{thm:policy_convergence}]
		For notation simplicity, set $\Delta^n_t:=\Delta^{\lambda^n}(t,X_t)$ and $P_t^{n}:=P_t^{\lambda^n}$ for $t\in[0,T]$, where 
		$(P^{\lambda^n},Q^{\lambda^n})=(P_t^{\lambda^n},Q_t^{\lambda^n})_{t\in[0,T]}$ is the solution of FBSDE \eqref{eq:BSDE} with $\lambda$ being replaced by $\lambda^n$. Then, it follows from Proposition~\ref{prop:performancediff} that
		the difference process $\delta P_t^n:=P_t^{n+1}-P_t^n$ satisfies that, for any $n\geq0$,
		\begin{align}\label{eq:policy_difference}
			\delta P^{n}_t=\E\left[\int_t^T\exp\left(-\int_t^s\lambda^{n+1}_r\d r\right)(\lambda^{n+1}_s-\lambda_s^n)\Delta^{n}_s\d s\Big|\F_t^1\right],\quad \forall t\in [0,T].
		\end{align}
		By using the policy improvement result given in Lemma~\ref{lem:policy_improvement}, we have $\delta P^n=(\delta P_t^n)_{t\in[0,T]}$ is a nonnegative process. This yields that $P_t^{n+1}\geq P_t^{n}$, and hence $\Delta^{n}_t=g(t,X_t)-P_t^{n}\leq g(t,X_t)-P_t^{n-1}=\Delta^{n-1}_t$ for any $n\geq1$.
		Thus, we have from the iterative form \eqref{eq:iterative000} that, for all $n\geq 0$,
		\begin{align}\label{eq:delta_B}
			(\lambda^{n+1}_t-\lambda^n_t)\Delta^{n}_t&=(M-M^{-1})\boldsymbol{1}_{\{\Delta^{n-1}_t\geq 0\}}(-\Delta_t^{n})^+=(M-M^{-1})\boldsymbol{1}_{\{\Delta^{n-1}_t\geq 0\}}\left(\delta P^{n-1}_t-\Delta^{n-1}_t\right)^+\nonumber\\
			&\leq (M-M^{-1})\boldsymbol{1}_{\{\Delta^{n-1}_t\geq 0\}}\delta P_t^{n-1}.
		\end{align}
		The last inequality follows from the fact that $(\delta P^{n-1}_t-\Delta^{n-1}_t)^+\leq \delta P_t^{n-1}$ on the event $\{\Delta^{n-1}_t\geq 0\}$. Plugging \eqref{eq:delta_B} into \eqref{eq:policy_difference}, we deduce that, for any $n\geq1$,
		\begin{align*}
			\delta P_t^{n}\leq M\E\left[\int_t^T\delta P_s^{n-1}\d s\Big|\F_t^1\right],\quad \forall t\in[0,T].
		\end{align*}
		Taking expectations on both sides of the above display, we have from Fubini's theorem that
		\begin{align*}
			\E\left[\delta P_t^n\right]\leq M\int_t^T\E\left[\delta P_s^{n-1}\right]\d s\leq M^n\int_t^T\frac{(s-t)^{n-1}}{(n-1)!}\E\left[\delta P_s^0\right]\d s.
		\end{align*}
		Set $C_0:=\int_0^T\E\left[\delta P_s^0\right]\d s<\infty$. Then, one has, for any $n\geq1$,
		\begin{align}\label{eq:policy_convergence}
			\mathcal{J}(\lambda^{n+1})-\mathcal{J}(\lambda^n)=\E\left[\delta P_0^n\right]\leq C_0M\frac{(MT)^{n-1}}{(n-1)!}
		\end{align}
		This yields that $\lim_{n\to\infty}\mathcal{J}(\lambda^n)=\mathcal{J}^*$ exists and the following convergence rate holds true, for some constant $C>0$ independent of $n$,
		\begin{align*}
			\left|\mathcal{J}(\lambda^n)-\mathcal{J}^*\right|\leq C\frac{(MT)^n}{n!}.
		\end{align*}
		Since $n\to\Delta^{n}_t$ is decreasing, the uniformly bounded intensity $n\to\lambda_t^n$ is decreasing, and hence admits a limit $\lambda_t^*$. Set $\Delta_t^*:=\Delta^{\lambda^*}(t,X_t)$ and $P_t^*=P_t^{\lambda^*}$ for $t\in[0,T]$. In light of Lemma~\ref{lem:delta_stability}, we have, as $n\to\infty$, 
		\begin{align*}
			\E\left[\sup_{t\in [0,T]}\left|\Delta^{*}_t-\Delta_t^{n}\right|^2\right]= \E\left[
			\sup_{t\in [0,T]}\left|P_t^{*}-P_t^{n}\right|^2\right]\leq C\|\lambda^n-\lambda^*\|_{L^2}^2\to 0,
		\end{align*}
		where, the last convergence follows from the bounded convergence theorem. Hence, there exists a subsequence $(\lambda^{n_k})_{k\geq1}\subset(\lambda^n)_{n\geq1}$ such that $\sup_{t\in [0,T]}|\Delta_t^{*}-\Delta_t^{n}|\to 0$, as $n\to\infty$, $\Pb$-a.s.. If $\Delta_t^{*}\geq0$, by the monotonicity of $n\to\Delta_t^{n}$, we have 
		$\Delta_t^{n_k}\geq0$ for every $k\geq1$, and thus $\lambda_t^{n_k+1}=M$. Passing to the limit gives $\lambda_t^*=M$.
		On the other hand, if $\Delta_t^{*}<0$, then
		$\Delta_t^{n_k}<0$ for all sufficiently large $k\geq1$,
		which implies $\lambda_t^{n_k+1}=M^{-1}$. 
		Consequently, we have $\lambda_t^*=M^{-1}$. This gives that $\lambda_t^*=M\mathbf 1_{\{\Delta_t^{*}\geq0\}}+\frac1M\mathbf 1_{\{\Delta_t^{*}<0\}}$, $\d t\times\d\mathbb{P}$-a.s.. This together with Proposition~\ref{prop:performancediff} yields that $\lambda^*=(\lambda_t^*)_{t\in[0,T]}\in\mathscr{U}^{M,{\rm m}}[0,T]$ is an optimal intensity. By using Lemma~\ref{lem:delta_stability}, we further have, as $n\to\infty$,
		\begin{align*}
			\left|\mathcal{J}(\lambda^n)-\mathcal{J}(\lambda^*)\right|^2\leq \E\left[
			\sup_{t\in [0,T]}\left|P_t^{*}-P_t^{n}\right|^2\right]\leq C\|\lambda^n-\lambda^*\|_{L^2}^2\to 0,
		\end{align*}
		where, the constant $C>0$ is independent of $n$.  This concludes that $\mathcal{J}^*=\mathcal{J}(\lambda^*)=\mathcal{V}^M$. Thus,  the proof of the theorem is complete.
	\end{proof}
	
	\begin{proof}[Proof of Theorem \ref{prop:relaxed-policy-convergence}]
		To ease the notation, we set $d_t^n:=M\mathbf 1_{\{h^{\theta^*(\lambda^n)}(t,X_t)\geq0\}}+M^{-1}\boldsymbol{1}_{\{h^{\theta^*(\lambda^n)}(t,X_t)<0\}}-\lambda_t^n$ with $\lambda_t^n:=\lambda^n(t,X_t)$ for $n\geq 0$. 
		Then, $\lambda_t^{n+1}=\lambda_t^n+\rho_nd_t^n$ and $|d_t^n|\leq M$. Set $\Delta_t^{n}:=\Delta^{\lambda^n}(t,X_t)$, $L_t^{n}:=L_t^{\lambda^n}$ and $e_t^n:=h^{\theta^*(\lambda^n)}(t,X_t)-\Delta_t^{n}$ for $t\in[0,T]$. Then, we have
		\begin{align*}
			d_t^n\Delta_t^{n}\geq(M-\lambda_t^n)(\Delta_t^{n})^++(\lambda_t^n-M^{-1})(\Delta_t^{n})^--M|e_t^n|.
		\end{align*}
		In fact, if the signs of $h^{\theta^*(\lambda^n)}(t,X_t)$
		and $\Delta_t^{n}$ agree, the first two terms on RHS of the above display give exactly $d_t^n\Delta_t^{n}$. If their signs differ, then
		$|\Delta_t^{n}|\leq |e_t^n|$, which yields the above inequality. Hence, we have from Proposition~\ref{prop:performancediff} that
		\begin{align}\label{eq:update_dif}
			\mathcal J(\lambda^{n+1})-\mathcal J(\lambda^n)=\rho_n\mathbb E\left[\int_0^T
			L_t^{n+1}d_t^n\Delta_t^{n}\,\mathrm dt\right].
		\end{align}
		Since $L_t^{n+1}/L_t^{n}=\exp(-\rho_n\int_0^t d_s^n\d s)$, one has $e^{-\rho_nMT}L_t^{n}
		\leq L_t^{n+1}\leq e^{\rho_nMT}L_t^{n}$. Combining the above equalities, we get
		\begin{align*}
			\mathcal J(\lambda^{n+1})-\mathcal J(\lambda^n)\geq
			\rho_ne^{-\rho_nMT}\mathcal R_n^M - M\rho_n\mathbb E\left[\int_0^T
			L_t^{n+1}|e_t^n|\d t\right].
		\end{align*}
		Hence, the Cauchy--Schwarz inequality implies that 
		\begin{align*}
			\mathbb E\left[\int_0^T
			L_t^{n+1}|e_t^n|\d t\right]&\leq
			\left\{\mathbb E\left[\int_0^T L_t^{n+1}\d t\right]\right\}^{1/2}\left\{\mathbb E\left[\int_0^TL_t^{n+1}|e_t^n|^2\d t\right]\right\}^{1/2}\nonumber\\
			&\leq\sqrt{2T} e^{\rho_nMT/2}\sqrt{S_M^{\lambda^n}(\theta^*)}.
		\end{align*}
		Consequently, one has
		\begin{equation}\label{eq:policy-improvement-estimate}
			\mathcal J(\lambda^{n+1})-\mathcal J(\lambda^n)\geq 
			\rho_ne^{-\rho_nMT}\mathcal R_n^M - M\sqrt{2T}\rho_ne^{\rho_nMT/2}
			\sqrt{S_M^{\lambda^n}(\theta^*)}.
		\end{equation}
		In particular, since $0<\rho_n\leq1$, it holds that
		\begin{align}\label{eq:uniform-policy-improvement}
			\mathcal J(\lambda^{n+1})-\mathcal J(\lambda^n)\geq
			e^{-MT}\rho_n\mathcal R_n^M -M\sqrt{2T}e^{MT/2}\rho_n\sqrt{S_M^{\lambda^n}(\theta^*)}.
		\end{align}
		Then, for any $N\geq1$, summing \eqref{eq:uniform-policy-improvement} from $n=0$ to $N-1$, we obtain
		\begin{align*}
			e^{-MT}\sum_{n=0}^{N-1}\rho_n\mathcal R_n^M\leq\mathcal J(\lambda^N)-\mathcal J(\lambda^0)+M\sqrt{2T}e^{MT/2}\sum_{n=0}^{N-1}\rho_n\sqrt{S_M^{\lambda^n}(\theta^*)}.
		\end{align*}
		Note that $\mathcal{J}(\lambda^N)\leq\mathcal{V}^M<\infty$ is upper bounded. Then, the learning error condition implies that $\sum_{n=0}^\infty\rho_n\mathcal R_n^M<\infty$.
		
		Let $\lambda,\mu\in\mathcal \U^{M,{\rm m}}[0,T]$. Applying Lemma \ref{lem:delta_stability} to solutions $P^{\lambda}$ and $P^{\mu}$ of the linear BSDE \eqref{eq:BSDE} under respective intensity processes $\lambda$ and $\mu$, there exists a constant $C_M>0$ independent of $(\lambda,\mu)$ such that
		\begin{align}\label{eq:delta-stability}
			\left\|\Delta^{\mu}-\Delta^{\lambda}\right\|_{L^2}\leq C_M\left\|\mu-\lambda\right\|_{L^2}.
		\end{align}
		Moreover, $|L_t^\mu-L_t^\lambda|
		\leq\int_0^t|\mu_s-\lambda_s|\d s$, which gives $\|L^\mu-L^\lambda\|_{L^2}
		\leq\frac{T}{\sqrt2}\|\mu-\lambda\|_{L^2}$. Define $b(a,x):=(M-a)x^++(a-M^{-1})x^-$ for $a\in[M^{-1},M]$ and $x\in\R$ and we have $
		|b(a,x)|\leq (M+M^{-1})|x|$ and $|b(a,x)-b(\tilde a,\tilde x)|\leq (M+M^{-1})|x-\tilde x|+|a-\tilde a||\tilde x|$ for $a,\tilde a\in [M^{-1},M]$ and $x,\tilde x\in\R$.
		This yields that
		\begin{align}\label{eq:R_dif}
			\left|\mathcal R^M(\mu)-\mathcal R^M(\lambda)\right|&\leq
			\mathbb E\left[\int_0^T|L_t^\mu-L_t^\lambda||b(\mu_t,\Delta_t^\mu)|\,\mathrm dt\right]+\mathbb E\left[\int_0^TL_t^\lambda
			|b(\mu_t,\Delta_t^\mu)-b(\lambda_t,\Delta_t^\lambda)|\d t\right]\nonumber\\
			&\quad\leq (M+M^{-1})\|L^\mu-L^\lambda\|_{L^2}
			\|\Delta^\mu\|_{L^2}+(M+M^{-1})\sqrt T
			\|\Delta^\mu-\Delta^\lambda\|_{L^2}\nonumber\\
			&\quad+\|\mu-\lambda\|_{L^2}
			\|\Delta^\lambda\|_{L^2}.
		\end{align}
		It then follows from Lemma~\ref{lem:delta_stability} that
		\begin{align}\label{eq:residual-lipschitz}
			\left|\mathcal R^M(\mu)-\mathcal R^M(\lambda)\right|\leq C_M\|\mu-\lambda\|_{L^2}.
		\end{align}
		In particular, since $\|\lambda^{n+1}-\lambda^n\|_{L^2}\leq M\sqrt T\rho_n$, 
		we have $|\mathcal R_{n+1}^M-\mathcal R_n^M|\leq C_M M\sqrt T\rho_n$. It then follows from Lemma A.1 in \cite{BoHuangWang26} that $\mathcal R_n^M\to 0$ as $n\to\infty$. 
		
		Taking the supremum over $\mu\in\mathcal \U^{M,{\rm m}}[0,T]$ in Proposition~\ref{prop:performancediff} yields that
		\begin{align}\label{eq:optimality-gap-residual}
			0\leq\mathcal{V}^M-\mathcal J(\lambda)=\E\left[\int_0^TL_t^{\lambda^{M,*}}(\lambda_t^{M,*}-\lambda_t)\Delta_t^{\lambda}\d t\right]
			\leq e^{MT}\mathcal R^M(\lambda).
		\end{align} 
		Here, the intensity process $\lambda^{M,*}=(\lambda_t^{M,*})_{t\in[0,T]}$ is given by \eqref{eq:lambdaMstardecoupling}. Taking $\lambda=\lambda^n$ in \eqref{eq:optimality-gap-residual} yields $0\leq\mathcal{V}^M-\mathcal J(\lambda^n)\leq e^{MT}\mathcal R_n^M\to0$ as $n\to\infty$ and the desired error decomposition \eqref{eq:error_decomp}. The proof of the proposition is complete.
	\end{proof}
	
	\begin{proof}[Proof of Lemma \ref{lem:learnQy}]
		It follows from It\^o's rule that
		\begin{align*}
			\d g(t,X_t)={\cal L}^{b,\sigma}g(t,X_t)\d t+\nabla_xg(t,X_t)^{\top}\sigma(t,X_t)\d W_t.   
		\end{align*}
		Then, for the reward process $\Phi=(\Phi_t)_{t\in[0,T]}$ defined in \eqref{eq:calQy}, one has, for $t\in[0,T]$,
		\begin{align}\label{eq:dPhitmart}
			\Phi_t=\Phi_0 + \int_0^t \mu(s,X_s) \d s + \int_0^t \nabla_xg(s,X_s)^{\top}\sigma(s,X_s)\d W_s.    
		\end{align}
		Set $A_t^\lambda:=\int_0^t\lambda(s,X_s)\d s$. We have from \eqref{eq:dPhitmart} and change of variables $z=A_s^\lambda$ (a localization argument is required if necessary) that 
		\begin{align*}
			{\cal Q}^{\lambda}(y)&=\E[\Phi_0]+\E\left[\int_0^{T}\mathbf{1}_{\{A_s^\lambda\leq y\}}\mu(s,X_s)\d s\right]=\E[\Phi_0]+\E\left[\int_0^{A_T^\lambda}\mathbf{1}_{\{z<y\}}\frac{\mu(\tau_z^\lambda,X_{\tau_z^\lambda})}{\lambda({\tau_z^\lambda},X_{\tau_z^\lambda})}\d z\right]\nonumber\\
			&=\E[\Phi_0]+\E\left[\int_0^{y}\mathbf{1}_{\{z<A_T^\lambda\}}\frac{\mu(\tau_z^\lambda,X_{\tau_z^\lambda})}{\lambda(\tau_z^\lambda,X_{\tau_z^\lambda})}\d z\right].
		\end{align*}
		This together with Fubini's theorem  yields the desired gradient formula.
	\end{proof}

\end{document}